\documentclass[final,leqno,a4paper]{amsart}

\usepackage[utf8]{inputenc}
\usepackage{amsmath}
\usepackage{amsfonts}
\usepackage{amssymb}
\usepackage{bbm}
\usepackage[font=footnotesize]{caption}
\usepackage[T1]{fontenc}
\usepackage{mathtools}
\usepackage{mathrsfs} 
\usepackage{amsbsy}

\newtheorem{remark}{Remark}

\usepackage{hyperref}

\usepackage{cleveref}
\usepackage{graphicx}
\usepackage{xcolor}
\usepackage[linesnumbered,ruled,vlined]{algorithm2e}
\newcommand{\SkipBeforeAndAfter}{\vspace{3mm}}

\SetCommentSty{mycommfont}

\usepackage{pgfplots}
\usepgfplotslibrary{groupplots}
\pgfplotsset{compat=1.18}

\def\letters{a,b,c,d,e,f,g,h,i,j,k,l,m,n,o,p,q,r,s,t,u,v,w,x,y,z,%
	A,B,C,D,E,F,G,H,I,J,K,L,M,N,O,P,Q,R,S,T,U,V,W,X,Y,Z}
\makeatletter
\@for \@l:=\letters \do{%
	\expandafter\edef\csname bb\@l\endcsname{\noexpand\ensuremath{\noexpand\mathbb{\@l}}}%
	\expandafter\edef\csname bf\@l\endcsname{{\noexpand\bf \@l}}%
	\expandafter\edef\csname cal\@l\endcsname{\noexpand\ensuremath{\noexpand\mathcal{\@l}}}%
	\expandafter\edef\csname eu\@l\endcsname{\noexpand\ensuremath{\noexpand\EuScript{\@l}}}%
	\expandafter\edef\csname frak\@l\endcsname{\noexpand\ensuremath{\noexpand\mathfrak{\@l}}}%
	\expandafter\edef\csname rm\@l\endcsname{{\noexpand\rm \@l}}%
	\expandafter\edef\csname scr\@l\endcsname{\noexpand\ensuremath{\noexpand\mathscr{\@l}}}%
}
\newcommand{\bfNull}{\boldsymbol{0}}
\newcommand{\bfSigma}{\boldsymbol{\Sigma}}
\newcommand{\bfLambda}{\boldsymbol{\Lambda}}
\newcommand{\bfscrE}{\pmb{\mathscr{E}}}
\newcommand{\bfscrL}{\pmb{\mathscr{L}}}
\newcommand{\bfscrS}{\pmb{\mathscr{S}}}

\newcommand{\dd}{\operatorname{d}\!}
\DeclareMathOperator*{\essinf}{ess\,inf}
\DeclareMathOperator*{\esssup}{ess\,sup}
\DeclareMathOperator\divergence{div}
\DeclareMathOperator{\rank}{rank}
\DeclareMathOperator{\diag}{diag}

\newcommand{\defis}{\mathrel{=\mathrel{\mathop:}}}

\newcommand{\Id}{\operatorname{Id}}
\renewcommand{\Re}{\operatorname{Re}}

\newcommand{\mytitle}{Multilevel Monte Carlo methods for random nonlinear eigenvalue problems}

\hypersetup{
    pdftitle={\mytitle},
    pdfauthor={J.~D\"olz and J.~Lozano-Sim\'on}
    }

\title[MLMC for random nonlinear eigenvalue problems]{\mytitle}

\author[J.~D\"olz and J.~Lozano-Sim\'on]{J.~D\"olz\protect\footnotemark[2]
        \and J.~Lozano-Sim\'on\protect\footnotemark[2]}
\makeatletter
\let\orig@makefnmark\@makefnmark
\g@addto@macro\@adminfootnotes{\let\@makefnmark\orig@makefnmark}
\makeatother

\thanks{Funded by the Deutsche Forschungsgemeinschaft (DFG, German Research Foundation) -- CRC 1720 -- 539309657 and  under Germany's Excellence Strategy -- EXC-2047/2 -- 390685813.}

\dedicatory{Dedicated to Professor Olaf Steinbach on the occasion of his 60th birthday.}

\begin{document}

\begin{abstract}
We develop a multilevel Monte Carlo (MLMC) estimator to compute statistical quantities of interest of solutions to parameter-dependent nonlinear eigenvalue problems. The main difficulty for this problem is that parame\-trized eigenpairs often lack a globally well-defined parameter-to-solution map due to crossings and bifurcations. These phenomena render eigenvalue-tracking infeasible in multi-dimensional parameter domains. To overcome this issue, we propose both an Eulerian and a Lagrangian perspective on parameter-to-solution maps for nonlinear eigenvalue problems, which remain well-defined even in the presence of local crossings and bifurcations. Their implementation is based on a straightforward modification of Beyn’s contour integral approach. Building on these maps, we propose an MLMC estimator that does not rely on matching eigenvalues or eigenvectors across samples or discretization levels. Nested approximation spaces and inter-level prolongation are only required when statistical quantities of interest related to the eigenspaces are of interest, but not for quantities of interest related to the eigenvalues. We demonstrate the feasibility and validity of our approach by numerical examples using a quasi-Monte Carlo method, including a nonlinear formulation of the Laplace Dirichlet eigenvalue problem on random domains in three dimensions motivated by the work of Steinbach and Unger. To the best of our knowledge, this is the first work considering the approximation of statistical quantities of interest for random nonlinear eigenvalue problems.
\end{abstract}

\maketitle

\footnotetext[2]{Institute for Numerical Simulation, University of Bonn, Germany.\\
\texttt{\{doelz,lozano\}@ins.uni-bonn.de}}

\section{Introduction}
\subsection{Objective}
We consider nonlinear eigenvalue problems as problems of finding a pair $(\lambda,u)\in\Lambda\times X$, with $\Lambda\subset\bbC$ and $X$ a suitable Hilbert space (such as $\bbC^n$), such that
\begin{equation}\label{eq:nonlinearevp}
\calA(\lambda) u = 0.
\end{equation}
Here, $\lambda\mapsto\calA(\lambda)$ is potentially nonlinear and $\calA(\lambda)\colon X \to X$ is a bounded linear operator (or, if $X=\bbC^n$, a matrix) for every value of $\lambda\in\Lambda$. The precise setting is detailed below. The problem \cref{eq:nonlinearevp} includes in particular linear eigenvalue problems
\begin{equation}\label{eq:linearEVP}
\calB u=\lambda u
\end{equation}
as a special case by setting $\calA(\lambda)=\calB-\lambda\Id$ and polynomial eigenvalue problems. Nonlinear eigenvalue problems arise in various applications such as structural mechanics, quantum physics, and acoustic problems, see \cite{GT2017,MV2004} for an overview.
A commonality of these examples is that their modelling is subject to random uncertainties due to measurement errors or manufacturing tolerances, for example. Accounting for these uncertainties leads to a random nonlinear eigenvalue problem of the form
\begin{equation}\label{eq:randomnonlinearevpintro}
\calA_\bfy(\lambda_\bfy) u_\bfy = 0,
\end{equation}
where $\bfy$ is a random parameter and $\calA_{\bfy}$ indicates dependence on $\bfy$. The solution of \cref{eq:randomnonlinearevpintro}, assuming it exists, is a random eigenpair $(\lambda_\bfy,u_\bfy)$.

We are interested in the numerical approximation of statistical quantities of interest of a few isolated eigenvalues located in a previously defined subset of $\Lambda$ and their corresponding eigenspaces. Exemplarily, we will consider suitable notions for the mean of the corresponding eigenvalues and eigenspaces, but our approach can be extended to many other statistical quantities of interest.

\subsection{Challenges}
The numerical approximation of said statistical quantities of interest is challenging for several reasons. First, on a theoretical level, it is already known for linear eigenvalue problems that parameter-to-solution mappings $\bfy\mapsto(\lambda_\bfy,u_\bfy)$ are not well defined in general, which makes the definition of any statistical quantity of interest problematic. This is due to the fact that the eigenvalues can cross, bifurcate as the parameter $\bfy$ varies, and the non-uniqueness of eigenvectors. Second, assuming that a well-defined meaning to these quantities can be given, the underlying integral is usually defined over a high-dimensional parameter domain, implying that its approximation using conventional quadrature rules suffers from the curse of dimensionality without further measures. Third, the (approximate) solution of \cref{eq:randomnonlinearevpintro} for a single parameter value $\bfy\in Y$ generally requires the solution of a large-scale matrix-valued nonlinear eigenvalue problem, which is computationally expensive. This makes any combination of conventional numerical methods for the approximation of statistical quantities of interest computationally prohibitive in most cases.

To the best of our knowledge, there is no existing work that addresses all of these challenges simultaneously for random nonlinear eigenvalue problems.

\subsection{Related work}
The analysis of parametric eigenvalue problems or eigenvalue problems with random input data has been mainly restricted to the linear eigenvalue problem \cref{eq:linearEVP} so far. Most of the early research has been done in structural analysis and aerospace engineering \cite{CT1969,SA1972b}, and we refer to the review papers \cite{AF2007,BR1988,PV1977} and the references therein for a detailed overview. Nevertheless, it seems that only eigenvalues of single multiplicity and the uncertainty quantification of rather small systems were discussed. Larger systems have only been made possible recently. For example, perturbation methods have been considered in \cite{DE2024,DESZ2024}, a sparse grid approach has been considered in \cite{AS2012}, stochastic collocation methods in \cite{GSH2023,HL2019a,GACS2019}, and (multilevel) quasi-Monte Carlo methods in \cite{CL2024,GGK+2019,GS2024a}. As an alternative to perturbation approaches and sampling-based methods, stochastic Galerkin methods were discussed in \cite{BOS2019,ES2019,GG2007,HL2019b,Wil2010}, with main emphasis on the acceleration of the eigensolvers. Let us also remark that parametrized eigenvalue problems have recently gained a lot of interest in the model order reduction community; see \cite{ABB2025,MSZ2026,TDE+2024} and the references therein.

Most of the aforementioned references restrict their attention to a single, isolated eigenvalue of single multiplicity and its corresponding eigenvector, which simplifies the definition of parameter-to-solution maps. For eigenvalues of higher multiplicity, or several eigenvalues, crossings and bifurcations can occur, which makes the definition of parameter-to-solution maps more difficult or even impossible already for linear eigenvalue problems \cite{Kat1995}. To overcome this issue, \cite{DE2024,DESZ2024,GSH2023,HL2019a} identify eigenspaces as the quantity of interest and derive local \cite{DE2024} or global \cite{GSH2023} analyticity results for eigenspaces. This allows to compute statistical quantities of interest by similar methods as single, isolated eigenvalues of single multiplicity or partial differential equations.

For parametric or stochastic nonlinear eigenvalue problems, the literature seems to be scarce. The only references considering parameter-to-solution maps for nonlinear eigenvalue problems we are aware of consider eigenvalue problems depending on a single, scalar-valued parameter. More precisely, in \cite{Xie2021}, a perturbation-based projection approach has been proposed to accelerate the computation of local parameter-to-invariant-pair maps for parameter-dependent nonlinear eigenvalue problems. A continuation approach based on invariant pairs to compute local parameter-to-solution maps for nonlinear eigenvalue problems was considered in \cite{BEK2011}. A more global perspective is taken in \cite{PB2024}, where an adaptive, non-intrusive, match-based collocation method is developed for approximating parameter-to-eigenvalue maps. The latter reference also includes strategies for dealing with eigenvalue migrations and bifurcations. In contrast to the setting considered here, these works do not address random nonlinear eigenvalue problems over high-dimensional parameter domains or the approximation of statistical quantities of eigenspaces.

\subsection{Contributions}
The goal of this work is to develop a computationally efficient algorithmic approach to approximate statistical quantities of interest of solutions to random nonlinear eigenvalue problems. The main contributions of this work are as follows:
\begin{enumerate}
	\item We propose an Eulerian and a Lagrangian perspective on parameter-to-solution maps for nonlinear eigenvalue problems. Both perspectives are well-defined even in the presence of local crossings and bifurcations. The implementation of those maps is a straightforward modification of Beyn's contour integral approach \cite{Bey2012}.
	\item We propose a multilevel Monte Carlo estimator for parametric nonlinear eigenvalue problems. Our proposed algorithm does not rely on any matching of eigenvalues or eigenvectors across samples or discretization levels. In fact, nestedness of approximation spaces and inter-level prolongation are only required if the mean of the eigenspace is of interest, but not for the mean of the eigenvalues. The extension to quasi-Monte Carlo and other quadrature methods is straightforward.
	\item We demonstrate the feasibility and validity of our approach on numerical examples based on a quasi-Monte Carlo method. Most prominently, we consider a nonlinear formulation of the Laplace Dirichlet eigenvalue problem on random domains in three dimensions. The latter example is motivated by the pioneering work of Steinbach and Unger, who proposed and analyzed the nonlinear formulation of the Laplace Dirichlet eigenvalue problem using boundary integral equations in \cite{SU2009,SU2012}.
\end{enumerate}
To the best of our knowledge, this is the first work considering the approximation of statistical quantities of interest for random nonlinear eigenvalue problems. The computational approaches for solving the arising nonlinear eigenvalue problems in \cite{SU2009,SU2012} have been the inspiration for the present work.

\subsection{Outline}
The remainder of this article is organized as follows. In \cref{sec:nonlinearevp}, we specify the setting of our considered nonlinear eigenvalue problems, their Galerkin discretization, Beyn's contour integral method, and our perspective on eigenspaces. \Cref{sec:p2s} introduces the Eulerian and Lagrangian parameter-to-solution maps and their discretization. In \cref{sec:MLMC}, we formulate the Monte Carlo and multilevel Monte Carlo estimators. Numerical experiments are presented in \cref{sec:experiments_toy} for a linear eigenvalue problem and in \cref{sec:experiments_BEM} for the nonlinear boundary integral formulation of the Laplace eigenvalue problem on random domains. \Cref{sec:conclusion} concludes the paper.

\section{Nonlinear eigenvalue problems}\label{sec:nonlinearevp}
\subsection{Eigenvalue problems for holomorphic Fredholm operator-valued functions}
In the following we recall a precise functional analytic setting for the nonlinear eigenvalue problem \cref{eq:nonlinearevp}. Let $\Lambda\subset\bbC$ be open and connected and let $X$ be a (complex) Hilbert space (such as $\bbC^n$). Denote by
\[
\calH^\infty\big(\overline{\Lambda};\calL(X)\big)
=
\{
f\colon\overline{\Lambda}\to \calL(X)\mid f~\text{is holomorphic and}~\|f\|_{\calH^\infty(\overline{\Lambda};\calL(X))}<\infty
\}
\]
the space of holomorphic, operator-valued functions on $\Lambda$ equipped with the supremum norm
\[
\|f\|_{\calH^\infty(\overline{\Lambda};\calL(X))}
=
\sup_{\lambda\in\overline{\Lambda}}\|f(\lambda)\|_{\calL(X)},
\]
with $\calL(X)$ denoting the space of bounded linear operators on $X$. Let
\[
\calA\in\calH^\infty(\overline{\Lambda};\calL(X))
\]
and assume that $\calA(\lambda)$ is a Fredholm operator of index zero for all $\lambda\in\Lambda$. We recall that a Fredholm operator of index zero is a bounded linear operator where the dimension of the kernel and codimension of the range coincide. All bounded, matrix-valued holomorphic functions satisfy these assumptions.
For the purpose of discretization, we will furthermore assume that the operator $\calA(\lambda)$ is coercive, see below.

The exposition of the following facts is closely aligned with \cite{SU2012}, see also \cite[Appendix]{KM1999}. The resolvent set of $\calA$ is given as
\[
\rho(\calA)
=
\{\lambda\in\Lambda\colon\text{there is}~\calA(\lambda)^{-1}\in\calL(X)\},
\]
which we assume to be non-empty. The spectrum of $\calA$ is given by $\sigma(\calA)=\Lambda\setminus\rho(\calA)$. For Fredholm operators of index zero, the spectrum consists of isolated points in $\Lambda$ and each $\lambda_0\in\Lambda$ is an eigenvalue of finite (geometric) multiplicity. More precisely, a number $\lambda_0\in\Lambda$ is called an eigenvalue of $\calA$ if there exists $u_0\in X\setminus\{0\}$ such that
\begin{equation}\label{eq:eigenpair}
\calA(\lambda_0)u_0=0.
\end{equation}
The (geometric) multiplicity of $\lambda_0$ is defined as the dimension of the kernel of $\calA(\lambda_0)$, denoted by $\ker(\calA(\lambda_0))$, which we call the eigenspace of $\lambda_0$. An ordered collection $u_0,\ldots,u_{m-1}\in X$ is called a Jordan chain of length $m$ for $\calA$ at $\lambda_0$ if
\begin{equation}\label{eq:Jordanchaindef}
\sum_{j=0}^n\frac{1}{j!}\calA^{(j)}(\lambda_0)u_{n-j}=0
\qquad
\text{for all}~n=0,\ldots,m-1,
\end{equation}
where $\calA^{(j)}(\lambda_0)$ denotes the $j$-th derivative of $\calA$ at $\lambda_0$. The maximal length of a Jordan chain to an eigenvalue $\lambda_0$ is denoted by $\varkappa(\calA,\lambda_0)$. The elements of any Jordan chain to an eigenvalue $\lambda_0$ are called generalized eigenelements of $\calA$ at $\lambda_0$ and their closed linear hull is called the generalized eigenspace of $\lambda_0$ and denoted by $G(\calA,\lambda_0)$. The dimension of the generalized eigenspace $G(\calA,\lambda_0)$ is called the algebraic multiplicity of $\lambda_0$. Eigenvalues with algebraic multiplicity one are called simple.

\subsection{Galerkin discretization}\label{sec:galerkin}
We recall the main results for the Galerkin discretization of nonlinear eigenvalue problems, following the exposition in \cite{SU2012}.
For the purpose of discretization we assume that $\calA(\lambda)$ is coercive, i.e., that it can be decomposed as
\[
\calA(\lambda) 
=
\calD+\calC(\lambda).
\]
Here, we assume that $\calD\in\calL(X)$ is $X$-elliptic, i.e., it satisfies
\[
\|x\|_X^2\lesssim\Re(\calD x,x)_X
\qquad
\text{for all}~x\in X,
\]
and that $\calC(\lambda)\in\calL(X)$ is compact for all $\lambda\in\Lambda$. It is easily verified that $\calA(\lambda)$ is a Fredholm operator of index zero for all $\lambda\in\Lambda$.

We consider a sequence of finite-dimensional subspaces $\{X_n\}_{n\in\bbN}\subset X$ with $N_n\coloneqq\dim X_n$ and $X_n\to X$ as $n\to\infty$ in the sense that
\[
\lim_{n\to\infty}\inf_{x_n\in X_n}\|x-x_n\|_X=0
\qquad
\text{for all}~x\in X.
\]
The variational eigenvalue problem for \cref{eq:eigenpair} is then to find $(\lambda_n,u_n)\in\Lambda\times (X_n\setminus\{0\})$ such that
\begin{equation}\label{eq:varform}
(\calA(\lambda_n)u_n,v_n)_X=0
\qquad
\text{for all}~v_n\in X_n.
\end{equation}

Denoting by $\bfu_n$ the coefficient vector of $u_n\in X_n$ with respect to a basis $\{\phi_1^{(n)},\ldots,\phi_{N_n}^{(n)}\}$ of $X_n$, \cref{eq:varform} can be rewritten as a finite-dimensional nonlinear eigenvalue problem
\begin{equation}\label{eq:discretevarform}
\bfA^{(n)}(\lambda_n)\bfu_n=\bfNull,
\end{equation}
with
\[
\lambda
\mapsto
\bfA^{(n)}(\lambda)\coloneqq\Big[\Big(\calA(\lambda)\phi_j^{(n)},\phi_i^{(n)}\Big)_X\Big]_{i,j=1,\ldots,N_n}
\in\bbC^{N_n\times N_n}
\]
being a matrix-valued holomorphic function. We drop the superscript ${}^{(n)}$ if it is not explicitly required from the context. By introducing the $X$-orthogonal projection $\calP_n\colon X\to X_n$, the continuous analogue to \cref{eq:discretevarform} is given by
\begin{equation}\label{eq:finitevarform}
	\calP_n\calA (\lambda_n)u_n=0.
\end{equation}

Based on the results of \cite{JW1977,Kar1996,Kar1996a} it was proven in \cite{SU2012} that for each eigenvalue $\lambda_0\in\sigma(\calA)$ there is a sequence of eigenvalues $\{\lambda_n\}_{n\in\bbN}$ of \cref{eq:discretevarform} such that $\lambda_n\to\lambda_0$ as $n\to\infty$. Moreover, for a sequence of eigenpairs $\{(\lambda_n,u_n)\}_{n\in\bbN}$ of \cref{eq:discretevarform} with $\|u_n\|_X=1$, $n\in\bbN$, there exists a subsequence of $\{u_n\}_{n\in\bbN}$ which converges to an eigenfunction $u_0$ of \cref{eq:eigenpair} with $\|u_0\|_X=1$. Spurious discrete eigenvalues can be excluded.

\subsection{Convergence rates}\label{sec:galerkinrates}
Let us recall how the approximation quality of the discrete eigenpair $(\lambda_n,u_n)$ to the continuous eigenpair $(\lambda_0,u_0)$ depends on the approximation capabilities of the finite-dimensional subspace $X_n$. To this end, we introduce the adjoint eigenvalue problem
\[
\calA^\star\big(\overline{\lambda}\big)v
\coloneqq
\big[\calA(\lambda)\big]^\star v=0
\]
and 
\begin{equation}\label{eq:delta}
\delta_n
\coloneqq
\max_{\substack{u_0\in G(\calA,\lambda_0)\\\|u_0\|_X\leq 1}}\inf_{u_n\in X_n}\|u_0-u_n\|_X,
\qquad
\delta_n^\star
\coloneqq
\max_{\substack{v_0\in G(\calA^\star,\overline{\lambda_0})\\\|v_0\|_X\leq 1}}\inf_{v_n\in X_n}\|v_0-v_n\|_X.
\end{equation}
Let moreover $\Lambda_c\subset\Lambda$ be compact with $\partial\Lambda_c\subset\rho(\calA)$ and $\Lambda_c\cap\sigma(\calA)=\{\lambda_0\}$. Then there is $N\in\bbN$ such that
\begin{equation}\label{eq:eigenvalueconvergence}
|\lambda_0-\lambda_n|
\lesssim
(\delta_n\delta_n^\star)^{1/\varkappa(\calA,\lambda_0)}
\end{equation}
for all $n\geq N$ and $\lambda_n\in\sigma(P_n\calA|_{X_n})\cap\Lambda_c$. The corresponding eigenelements $u_n$ normalized to $\|u_n\|_X=1$ satisfy
\begin{equation}\label{eq:eigenelementconvergence}
\inf_{u_0\in\ker(\calA(\lambda_0))}\|u_0-u_n\|_X 
\lesssim
|\lambda_0-\lambda_n|
+
\max_{\substack{x_0\in\ker(\calA(\lambda_0))\\\|x_0\|_X\leq 1}}\inf_{x_n\in X_n}\|x_0-x_n\|_X.
\end{equation}
The hidden constants in \cref{eq:eigenvalueconvergence,eq:eigenelementconvergence} only depend on the multiplicities of the eigenvalue $\lambda_0$, the contour $\partial\Lambda_c$, and $\operatorname{dist}(\lambda_0,\partial\Lambda_c)$, $\|\calA\|_{\calH^\infty(\overline{\Lambda_c},\calL(X))}$, and $\|\calA^{-1}\|_{\calH^\infty(\partial\Lambda_c,\calL(X))}$. We refer to \cite{SU2012,Ung2009} for further details.

\subsection{Beyn's contour integral approach}\label{sec:Beynoriginal}
For the numerical solution of the discretized eigenvalue problem \cref{eq:discretevarform} we use the contour integral approach from Beyn \cite{Bey2012}. It aims at finding the generalized eigenspaces to a small number of eigenpairs whose eigenvalues are contained within a predefined smooth contour $\Gamma\subset\rho(P_n\calA|_{X_n})$ enclosing the bounded set $\Lambda_c$. The underlying principle is to use Keldysh's theorem to reduce the nonlinear eigenvalue problem \cref{eq:discretevarform} to a linear $m\times m$-matrix eigenvalue problem, where
\[
m
=
\sum_{\lambda_0\in\sigma(P_n\calA|_{X_n})\cap\Lambda_c}\dim(G(P_n\calA|_{X_n},\lambda_0))
\]
is the number of eigenvalues enclosed by $\Gamma$, counting with algebraic multiplicity. In practice, $m$ can be suitably estimated during the algorithm and only the contour $\Gamma$ needs to be provided. We refer to \cite{Bey2012} for details and convergence estimates and only recall the main algorithm.

The starting point of the approach is a matrix $\widehat{\bfV}\in\bbC^{N_n\times m'}$ of full rank and $m'\geq m$. The approach further considers the integrals
\begin{align*}
\bfT_0
&\coloneqq\calP_{G(\bfA,\Lambda_c)}\widehat{\bfV}
\coloneqq
\frac{1}{2\pi i}\int_{\Gamma}\bfA(z)^{-1}\widehat{\bfV}\dd z\in\bbC^{N_n\times m'},\\
\bfT_1
&\coloneqq\phantom{\coloneqq\calP_{G(\bfA,\Lambda_c)}\widehat{\bfV}}
\frac{1}{2\pi i}\int_{\Gamma}z\bfA(z)^{-1}\widehat{\bfV}\dd z\in\bbC^{N_n\times m'},
\end{align*}
assuming that $\rank(\bfT_0)=m$. In practice, this matrix can be reliably achieved by choosing $\widehat{\bfV}$ as a random matrix. Moreover, in algorithmic implementations of the approach, these integrals are approximated using the trapezoidal rule. With $\Gamma$ being smooth, the trapezoidal rule converges exponentially, and only a few quadrature points are required \cite{Bey2012}. We neglect the quadrature error in the following. For later reference we also note that the columns of $\bfT_0=\calP_{G(\bfA,\Lambda_c)}\widehat{\bfV}$ lie in the direct sum of the generalized eigenspaces of the eigenvalues in $\Lambda_c$, i.e., in
\[
G(\bfA,\Lambda_c)
=
\bigoplus_{\lambda_0\in\sigma(\bfA)\cap\Lambda_c}G(\bfA,\lambda_0).
\]
For linear eigenvalue problems, $\calP_{G(\bfA,\Lambda_c)}$ coincides with the Riesz projection.
One can show that $\bfT_0$ has (numerical) rank $m$ and, using Keldysh's theorem, that
\begin{equation}\label{eq:T0T1}
\bfT_0=\bfV\bfW^\star\widehat{\bfV},
\qquad
\bfT_1=\bfV\bfLambda\bfW^\star\widehat{\bfV},
\end{equation}
where $\bfV,\bfW\in\bbC^{N_n\times m}$ and $\bfLambda$ is a Jordan matrix to the eigenvalues in $\sigma(\bfA)\cap\Lambda_c)$, including multiplicities. Note further that
\begin{equation}\label{eq:GspanT}
G(\bfA,\Lambda_c)
=
\operatorname{span}(\bfT_0)
=
\operatorname{span}(\bfV).
\end{equation}
The goal of Beyn's contour integral approach is to recover $\bfLambda$ and $\bfV$ from $\bfT_0$ and $\bfT_1$ using algebraic manipulations to obtain the eigenvalues in $\Lambda_c$, their multiplicities, and a basis of the corresponding generalized eigenspace.

The traditional form of Beyn's algorithm is to let $\bfT_0=\bfV_0\bfSigma_0\bfW_0^\star$ the reduced SVD of $\bfT_0$ and to observe that $\operatorname{span}(\bfV)=\operatorname{span}(\bfV_0)$ implies that there is a nonsingular matrix $\bfS$ such that $\bfV=\bfV_0\bfS$ and thus $\bfV_0^\star\bfV=\bfS$. Hence, $\bfW^\star\widehat{\bfV}=\bfS^{-1}\bfSigma_0\bfW_0^\star$ and $\bfV_0^\star\bfT_1=\bfS\bfLambda\bfW^\star\widehat{\bfV}$, implying that
\begin{equation}\label{eq:Beynevp}
	\bfB
	=
	\bfS\bfLambda\bfS^{-1}
	=
	\bfV_0^\star\bfT_1\bfW_0\bfSigma_0^{-1}\in\bbC^{m\times m}
\end{equation}
has the same eigenvalues as $\bfA$ has eigenvalues in $\Gamma$, including algebraic multiplicities. Further, we note that if $(\lambda_n,\bfs_n)$ is an eigenpair of $\bfB$, then $(\lambda_n,\bfV_0\bfs_n)$ is an eigenpair of \cref{eq:discretevarform}. Moreover, if $\bfs_n$ is a generalized eigenelement of $\bfB$ to the eigenvalue $\lambda_n$, then $\bfV_0\bfs_n$ is a generalized eigenelement of \cref{eq:discretevarform} to the eigenvalue $\lambda_n$. Thus, the approach can recover all eigenpairs and generalized eigenelements of \cref{eq:discretevarform} within the specified contour $\Gamma$. Note that we make no statements about orthogonality. The size of the linear eigenvalue problem \cref{eq:Beynevp} is independent of the chosen Galerkin discretization, and only depends on the geometry of $\sigma(\bfA)$ and $\Gamma$.

\subsection{A modified version of Beyn's approach}\label{sec:Beynmodified}
For the purpose of parameter dependent nonlinear eigenvalue problems such as \cref{eq:randomnonlinearevpintro}, an approach using the singular value decomposition is difficult to analyse. Thus, we will in the following provide an alternate approach for recovering $\bfLambda$ and a basis of the generalized eigenspaces $G(\bfA,\Lambda_c)$ from $\bfT_0$ and $\bfT_1$ in \cref{eq:T0T1}. The possibility of such an approach was already noted in Beyn's original article \cite{Bey2012}.

Observing that \cref{eq:GspanT} implies that there is a matrix $\bfR$ of rank $m$ (nonsingular for $m=m'$) such that $\bfV=\bfT_0\bfR$ yields that $\bfR=\bfT_0^\dagger\bfV$. Here, $\bfT_0^\dagger$ denotes the (unique) Moore-Penrose pseudoinverse of $\bfT_0$, see \cite{GV2013} for a definition. Hence, $\bfT_0=\bfT_0\bfR\bfW^\star\widehat{\bfV}$ and multiplying with $\bfT_0^\dagger$ implies $\bfW^\star\widehat{\bfV}=\bfR^{\dagger}$. This yields that
\begin{equation}\label{eq:Beynevp2}
\bfC
=
\bfT_0^\dagger\bfT_1
=
\bfT_0^\dagger\bfV\bfLambda\bfW^\star\widehat{\bfV}
=
\bfR\bfLambda\bfR^{\dagger}
\end{equation}
has the same nonzero eigenvalues as $\bfA$ in $\Lambda_c$, including algebraic multiplicities, and $m'-m$ additional zero eigenvalues. As before, we note that if $(\lambda_n,\bfr_n)$ is an eigenpair of $\bfC$, then $(\lambda_n,\bfT_0\bfr_n)$ is an eigenpair of \cref{eq:discretevarform}. The statement concerning generalized eigenelements holds verbatim. Moreover, the size of the linear eigenvalue problem \cref{eq:Beynevp2} is independent of the chosen Galerkin discretization, and only depends on the geometry of $\sigma(\bfA)$ and $\Gamma$. This observation is the key ingredient for our Monte Carlo and multilevel Monte Carlo approach.

\subsection{Nonlinear eigenvalue problems and eigenspaces}\label{sec:NEVPeigenspaces}
It is well known that linear eigenvalue problems \cref{eq:linearEVP} can straightforwardly be extended to eigenspaces. To this end, let $\{(\lambda_i,u_i)\}_{i=0}^{m-1}$ be eigenpairs to \cref{eq:linearEVP}. Setting $\bfu=[u_0,\ldots,u_{m-1}]$ and $\bfLambda=\diag(\lambda_0,\ldots,\lambda_{m-1})$ it holds
\[
\calB\bfu
\coloneqq
[\calB u_0,\ldots,\calB u_{m-1}]
=
\bfu\bfLambda.
\]
For any unitary matrix $\bfQ\in\bbC^{m\times m}$ and setting $\widetilde{\bfu}=\bfu\bfQ$ and $\widetilde{\bfLambda}=\bfQ^\star\bfLambda\bfQ$, this is equivalent to
\[
\calB\widetilde{\bfu}
=
\calB\bfu\bfQ
=
\bfu\bfQ\bfQ^\star\bfLambda\bfQ
=
\widetilde{\bfu}\widetilde{\bfLambda},
\]
with $\widetilde{\bfLambda}$ not necessarily diagonal.

In the following, we extend this interpretation to nonlinear eigenvalue problems where the eigenvalue dependence is holomorphic. The latter allows us to use a Taylor expansion of $\calA(\lambda)$ around a reference point $\lambda^\star$ and to rewrite \cref{eq:nonlinearevp} as
\[
0
=
\calA(\lambda)u
=
\sum_{\ell=0}^\infty\frac{1}{\ell!}\calA^{(\ell)}(\lambda^\star)u(\lambda-\lambda^\star)^\ell
\]
for all $\lambda$ such that $|\lambda-\lambda^\star|$ is smaller than the Taylor series' radius of convergence.
Setting $\bfu=[u_0,\ldots,u_{m-1}]$ and $\bfLambda=\diag(\lambda_0,\ldots,\lambda_{m-1})$ it holds
\[
\bfNull
=
[\calA(\lambda_0)u_0,\ldots,\calA(\lambda_{m-1})u_{m-1}]
=
\sum_{\ell=0}^\infty\frac{1}{\ell!}\calA^{(\ell)}(\lambda^\star)\bfu(\bfLambda-\lambda^\star\bfI)^\ell
\]
for all $\lambda_i$ such that $|\lambda_i-\lambda^\star|$ is smaller than the Taylor series' radius of convergence.
Then, for any unitary matrix $\bfQ\in\bbC^{m\times m}$ and setting $\widetilde{\bfu}=\bfu\bfQ$ and $\widetilde{\bfLambda}=\bfQ^\star\bfLambda\bfQ$, this is equivalent to
\begin{equation}\label{eq:nonlineareigenmatrixproblem}
\bfNull
=
\sum_{\ell=0}^\infty\frac{1}{\ell!}\calA^{(\ell)}(\lambda^\star)\widetilde{\bfu}(\widetilde{\bfLambda}-\lambda^\star\bfI)^\ell
\eqqcolon
\calA(\widetilde{\bfLambda})\widetilde{\bfu}.
\end{equation}
This extends the nonlinear eigenvalue problem to the setting of eigenspaces. We call $\widetilde{\bfLambda}$ an eigenmatrix to the eigenspace spanned by $\widetilde{\bfu}$.

To generalize this interpretation from eigenspaces to generalized eigenspaces, consider $k$ eigenvalues $\lambda_0,\ldots,\lambda_{k-1}$ of $\calA$, counting with geometric multiplicities. Denote the corresponding algebraic multiplicities by $\varkappa_0,\ldots,\varkappa_{k-1}$, assuming that the corresponding generalized eigenspaces are linearly independent, and denote by $m=\sum_{i=0}^{k-1}\varkappa_i$ the total algebraic multiplicity of all eigenvalues in $\Lambda_c$. Write $\bfu=[u_0,\ldots,u_{m-1}]$ and define the corresponding Jordan matrix as
\[
\bfJ = \begin{bmatrix}
\bfJ_{\varkappa_0}(\lambda_0) & & 0 \\
& \ddots & \\
0 & & \bfJ_{\varkappa_{k-1}}(\lambda_{k-1})
\end{bmatrix},
\qquad
\bfJ_{\varkappa}(\lambda) = \begin{bmatrix}
\lambda & 1 & & 0 \\
& \lambda & \ddots & \\
& & \ddots & 1 \\
0 & & & \lambda
\end{bmatrix} \in \mathbb{C}^{\varkappa \times \varkappa},
\]
The nonlinear eigenmatrix problem then becomes
\begin{equation}\label{eq:nonlineareigenmatrixproblem_jordan}
\sum_{\ell=0}^\infty\frac{1}{\ell!}\calA^{(\ell)}(\lambda^\star)\bfu(\bfJ-\lambda^\star\bfI)^\ell = \bfNull,
\end{equation}
which is the natural generalization of \cref{eq:nonlineareigenmatrixproblem} to the case of higher algebraic multiplicities whenever the Taylor series converges. Note that in case of a single eigenvalue $\lambda_0=\lambda^\star$ recovers the defining property of a Jordan chain from \cref{eq:Jordanchaindef}.

Several remarks to \cref{eq:nonlineareigenmatrixproblem} are in order. First, we note that our discussions above reveal that for $m>1$ solutions to \cref{eq:nonlineareigenmatrixproblem} are never unique. In fact, $(\widetilde{\bfLambda},\widetilde{\bfu})$ could contain $m$ copies of the same eigenpair, such that further conditions are required to obtain all solutions. Second, our discussion above also reveals that the eigenmatrix $\widetilde{\bfLambda}$ of solutions to \cref{eq:nonlineareigenmatrixproblem} is not necessarily diagonal. Finally, assuming that $(\widetilde{\bfLambda},\widetilde{\bfu})$ is a solution to \cref{eq:nonlineareigenmatrixproblem} and $\widetilde{\bfLambda}=\bfR\bfLambda\bfR^{-1}$, $\bfLambda=\diag(\lambda_0,\ldots,\lambda_{m-1})$, for some nonsingular matrix $\bfR\in\bbC^{m\times m}$, then $(\bfLambda,\widetilde{\bfu}\bfR)=(\bfLambda,\bfu)$ with $\bfu=[u_0,\ldots,u_{m-1}]$ yields solutions $(\lambda_i,u_i)$, $i=0,\ldots,m-1$, to \cref{eq:nonlinearevp}.

We will exploit these observations later to define parameter-to-solution maps for random nonlinear eigenvalue problems.

\section{Parameter-to-solution maps for random nonlinear eigenvalue problems}\label{sec:p2s}
\subsection{Random nonlinear eigenvalue problems}\label{sec:randomnonlinearevp}
Let $(Y,\calF,\bbP)$ be a probability space and let $\Lambda\subset\bbC$ be open and connected. We consider $\calA$ from \cref{eq:eigenpair} to be dependent on a random parameter $\bfy\in Y$ and, for simplicity, assume that the random variable $\bfy\mapsto\calA_\bfy$ is essentially bounded, i.e., that
\[
\calA_{\bullet}\in L_\bbP^\infty\Big(Y;\calH^\infty\big(\overline{\Lambda};\calL(X)\big)\Big).
\]
We consider the parametrized nonlinear eigenvalue problem
\begin{equation}\label{eq:randomnonlinearevp}
\calA_{\bfy}(\lambda_\bfy)u_\bfy=0,
\qquad
\bfy\in Y,
\end{equation}
and readily remark that all considerations from \cref{sec:nonlinearevp} apply pointwise to \cref{eq:randomnonlinearevp} for $\bbP$-almost all $\bfy\in Y$. Let $\Lambda_c\subset\Lambda$ be open and bounded and set $\Gamma=\partial\Lambda_c$. We assume the existence of a small number $m\in\bbN$ of eigenvalues within $\Lambda_c$ for $\bbP$-almost all $\bfy\in Y$, counting with algebraic multiplicity. We assume all other eigenvalues to be in $\Lambda_c^c=\Lambda\setminus\overline{\Lambda_c}$ and that
\begin{align*}
0&<\essinf_{\bfy\in Y}\operatorname{dist}\big(\sigma(\calA_{\bfy})\cap\Lambda_c,\Gamma\big)
=
\essinf_{\bfy\in Y}\operatorname{dist}\big(\sigma(\calA_{\bfy})\cap\Lambda_c,\Lambda_c^c\big),\\
0&<\essinf_{\bfy\in Y}\operatorname{dist}\big(\sigma(\calA_{\bfy})\cap\Lambda_c^c,\Gamma\big)=\essinf_{\bfy\in Y}\operatorname{dist}\big(\sigma(\calA_{\bfy})\cap\Lambda_c^c,\Lambda_c\big),
\end{align*}
i.e., we impose a uniform minimal distance between all eigenvalues and $\Gamma$.  Moreover, we assume that
\[
\calA_{\bullet}(\bullet)^{-1}\in L_\bbP^\infty\Big(Y;\calH^\infty\big(\Gamma;\calL(X)\big)\Big),
\]
i.e., that the resolvent is uniformly bounded for all $\lambda\in\Gamma$ and $\bbP$-almost all $\bfy\in Y$, and that the threshold $N$ and the hidden constants in \cref{eq:eigenvalueconvergence,eq:eigenelementconvergence} can be chosen independently of $\bfy$ for $\bbP$-almost all $\bfy\in Y$.

Now, consider eigenpairs $(\lambda_{0,\bfy},u_{0,\bfy}),\ldots,(\lambda_{m-1,\bfy},u_{m-1,\bfy})\subset\Lambda_c\times X$ of \cref{eq:randomnonlinearevp} for $\bbP$-almost all $\bfy\in Y$. These eigenpairs are not unique, as the eigenvalues $\lambda_{0,\bfy},\ldots,\lambda_{m-1,\bfy}$ can be ordered arbitrarily and the eigenfunctions $u_{0,\bfy},\ldots,u_{m-1,\bfy}$ can be scaled arbitrarily. This and potential bifurcations of eigenvalues and eigenfunctions with respect to the parameter $\bfy$ make it difficult to define a parameter-to-solution map. To mitigate this issue, we observe that the generalized eigenspace $G(\calA_\bfy,\Lambda_c)$ is unique, but needs a basis representation for computational purposes. Fixing a basis representation for the generalized eigenspace $G(\calA_\bfy,\Lambda_c)$ for $\bbP$-almost all $\bfy\in Y$ is equivalent to defining a parameter-to-eigenspace map and readily enables us to define a parameter-to-eigenmatrix map as well. Combining these two maps yields a parameter-to-solution map for \cref{eq:randomnonlinearevp}.

In the following, we pursue two approaches, one Eulerian and one Lagrangian, to define parameter-to-solution maps for \cref{eq:randomnonlinearevp}. The approaches we pursue will also allow us to compute the parameter-to-eigenmatrix map without computing the parameter-to-eigenspace map explicitly.

\subsection{Eulerian parameter-to-solution maps}\label{sec:parameter-to-eigenspace-V1}

\subsubsection{Continuous perspective}
Our first approach for defining a parameter-to-solu\-tion map is based on an Eulerian perspective. That is, we fix a reference basis for the eigenspace, project the eigenvalue problem onto this fixed basis, and consider the eigenvalues and eigenmatrix with respect to this fixed basis. To fix a reference basis, we select a suitable reference value $\widehat{\bfy}\in Y$ and consider the eigenpairs $(\lambda_{0,\widehat{\bfy}},u_{0,\widehat{\bfy}}),\ldots,(\lambda_{m-1,\widehat{\bfy}},u_{m-1,\widehat{\bfy}})$ of \cref{eq:randomnonlinearevp} for $\bfy=\widehat{\bfy}$ with eigenvalues in $\Lambda_c$. This solves the eigenvalue problem
\begin{equation}\label{eq:referenceevp}
\calA_{\widehat{\bfy}}(\bfLambda_{0,\widehat{\bfy}})\bfu_{0,\widehat{\bfy}}=\boldsymbol{0},
\end{equation}
where $\bfu_{0,\widehat{\bfy}}=[u_{0,\widehat{\bfy}},\ldots,u_{m-1,\widehat{\bfy}}]$ and $\bfLambda_{0,\widehat{\bfy}}=\diag(\lambda_{0,\widehat{\bfy}},\ldots,\lambda_{m-1,\widehat{\bfy}})$. Projecting the eigenvalue problem \cref{eq:randomnonlinearevp} for $\bbP$-almost all $\bfy\in Y$ onto the reference basis $\bfu_{0,\widehat{\bfy}}$ corresponds to a pointwise Galerkin discretization, cf.\ \cref{sec:galerkin}, of \cref{eq:randomnonlinearevp} onto the finite dimensional subspace $\operatorname{span}(\bfu_{0,\widehat{\bfy}})\subset X$ such that \cref{eq:randomnonlinearevp} becomes
\begin{equation}\label{eq:projectednevp}
\calP_{\operatorname{span}(\bfu_{0,\widehat{\bfy}})}\calA_\bfy(\bfLambda_{\bfy,\widehat{\bfy}})\bfu_{\bfy,\widehat{\bfy}}=\boldsymbol{0},
\qquad
\operatorname{span}\{\bfu_{\bfy,\widehat{\bfy}}\}=\operatorname{span}\{\calP_{\operatorname{span}(\bfu_{0,\widehat{\bfy}})}\bfu_\bfy\}.
\end{equation}
We recall from \cref{sec:NEVPeigenspaces} that, for $m>1$, solutions to this eigenvalue problem are not unique anymore. To alleviate this issue, we observe that the condition on the spans of $\bfu_{\bfy,\widehat{\bfy}}$ and $\bfu_{0,\widehat{\bfy}}$ implies that there is a matrix $\bfR_\bfy$ such that $\bfu_{\bfy,\widehat{\bfy}}=\bfu_{0,\widehat{\bfy}}\bfR_\bfy$. Since we are only interested in a basis of the eigenspace, we may reformulate \cref{eq:projectednevp} to find $\bfB_\bfy=\bfR_\bfy\bfLambda_{\bfy,\widehat{\bfy}}\bfR_\bfy^{\dagger}\in\bbC^{m\times m}$ such that
\begin{equation}\label{eq:eulerianevp}
\calP_{\operatorname{span}(\bfu_{0,\widehat{\bfy}})}\calA_\bfy(\bfB_\bfy)\bfu_{0,\widehat{\bfy}}=\boldsymbol{0}.
\end{equation}
This resolves any issues concerning scaling and ordering of the eigenpairs and allows us to define a parameter-to-solution map for \cref{eq:randomnonlinearevp} as
\[
\scrE\colon Y\to\bbC^{m\times m}\times X^m,
\qquad
\bfy\mapsto\scrE(\bfy)\coloneqq(\bfB_\bfy,\bfu_{0,\widehat{\bfy}}).
\]
and the assumptions from \cref{sec:randomnonlinearevp} guarantee that
\[
\scrE\in L_\bbP^\infty\Big(Y;\bbC^{m\times m}\times X^m\Big),
\]
whenever the $m$-th singular value of $\bfV_{0,\widehat{\bfy}}^\star\bfT_{0,\bfy}$ is bounded uniformly from below for $\bbP$-almost all $y\in Y$.

\subsubsection{Discretization and implementation}
For the computational implementation of the Eulerian parameter-to-solution map, we consider a Galerkin discretization of \cref{eq:randomnonlinearevp} as in \cref{sec:galerkin}. Following Beyn's original contour integral approach, cf.\ \cref{sec:Beynoriginal}, we fix a matrix $\widehat{\bfV}\in\bbC^{N_n\times m'}$, $m'\geq m$, of full rank and with $\rank(\calP_{G(\bfA_\bfy,\Lambda_c)}\widehat{\bfV})=m$ for $\bbP$-almost all $\bfy\in Y$, and introduce the parameter dependent integrals
\begin{subequations}\label{eq:T01y}
\begin{align}
	\bfT_{0,\bfy}
	&\coloneqq\calP_{G(\bfA_\bfy,\Lambda_c)}\widehat{\bfV}\coloneqq
	\frac{1}{2\pi i}\int_{\Gamma}\bfA_\bfy(z)^{-1}\widehat{\bfV}\dd z\in\bbC^{N_n\times m'},\\
	\bfT_{1,\bfy}
	&\coloneqq\phantom{\coloneqq\calP_{G(\bfA_\bfy,\Lambda_c)}\widehat{\bfV}}
	\frac{1}{2\pi i}\int_{\Gamma}z\bfA_\bfy(z)^{-1}\widehat{\bfV}\dd z\in\bbC^{N_n\times m'},
\end{align}
\end{subequations}
for all $\bfy\in Y$. 

Now, for a fixed reference value $\widehat{\bfy}\in Y$, the discrete analog to \cref{eq:referenceevp} is given by
\[
\bfA_{\widehat{\bfy}}(\bfLambda_{0,\widehat{\bfy}})\bfV_{0,\widehat{\bfy}}=\bfNull,
\]
and $\bfV_{0,\widehat{\bfy}}$ can be obtained from the reduced SVD $\bfT_{0,\widehat{\bfy}}=\bfV_{0,\widehat{\bfy}}\bfSigma_{0,\widehat{\bfy}}\bfW_{0,\widehat{\bfy}}^\star$ of $\bfT_{0,\widehat{\bfy}}$. Having $\bfV_{0,\widehat{\bfy}}$ at hand, we can project the discrete eigenvalue problem \cref{eq:discretevarform} for $\bbP$-almost all $\bfy\in Y$ onto the reference basis $\bfV_{0,\widehat{\bfy}}$ to obtain a discrete version of \cref{eq:eulerianevp} as
\[
\bfV_{0,\widehat{\bfy}}\bfV_{0,\widehat{\bfy}}^\star\bfA_\bfy(\bfLambda_\bfy)\bfu_{0,\bfy}=\bfNull,
\qquad
\operatorname{span}\{\bfu_{0,\bfy}\}=\operatorname{span}\{\bfV_{0,\widehat{\bfy}}\bfV_{0,\widehat{\bfy}}^\star\bfV_\bfy\}.
\]
Now, letting $\bfT_{0,\bfy}=\bfV_{0,\bfy}\bfSigma_{0,\bfy}\bfW_{0,\bfy}^\star$ be the reduced SVD of $\bfT_{0,\bfy}$ and observing that $\bfV_{0,\widehat{\bfy}}$ is an orthonormal basis of
\[
G(\bfV_{0,\widehat{\bfy}}\bfV_{0,\widehat{\bfy}}^\star\bfA,\Lambda_c)
=
\operatorname{span}(\bfV_{0,\widehat{\bfy}})
=
\operatorname{span}(\bfV_{0,\widehat{\bfy}}\bfV_{0,\widehat{\bfy}}^\star\bfV_{0,\bfy})
=
\operatorname{span}(\bfV_{0,\widehat{\bfy}}\bfV_{0,\widehat{\bfy}}^\star\bfT_{0,\bfy}),
\]
implies that there is a matrix $\bfR_\bfy\in\bbC^{m\times m}$ such that $\bfV_{0,\widehat{\bfy}}\bfV_{0,\widehat{\bfy}}^\star\bfV_{0,\bfy}=\bfV_{0,\widehat{\bfy}}\bfR_\bfy$. In complete analogy to the derivation in \cref{sec:Beynoriginal}, we derive
\begin{align*}
\bfSigma_{0,\bfy}&=\bfV_{0,\widehat{\bfy}}^\star\bfT_{0,\bfy},\\
\bfB_\bfy&=\bfV_{0,\widehat{\bfy}}^\star\bfT_{1,\bfy}\bfSigma_{0,\bfy}^{\dagger}.
\end{align*}
If $\bfSigma_{0,\bfy}$ has full rank $m=m'$, then $\bfSigma_{0,\bfy}^{\dagger}=\bfSigma_{0,\bfy}^{-1}$ and the derivation implies that $\bfB_\bfy$ is a similarity transformation of \cref{eq:Beynevp}, i.e., the matrix $\bfB$ when the contour integral method is applied to \cref{eq:referenceevp} at $\bfy$. This implies that $\bfB_\bfy$ has the same eigenvalues as $\bfA_\bfy$ has eigenvalues in $\Lambda_c$, including algebraic multiplicities, whenever $\rank(\bfSigma_{0.\bfy})=m=m'$.

Our discrete Eulerian parameter-to-solution map is thus given as
\[
\bfscrE^{(n)}\colon Y\to\bbC^{m\times m}\times\bbC^{N_n\times m},
\qquad
\bfy\mapsto\bfscrE^{(n)}(\bfy)\coloneqq(\bfB_\bfy,\bfV_{0,\widehat{\bfy}}),
\]
and the assumptions from \cref{sec:randomnonlinearevp} guarantee that
\[
\bfscrE^{(n)}\in L_\bbP^\infty\Big(Y;\bbC^{m\times m}\times\bbC^{N_n\times m}\Big),
\]
whenever the $m$-th singular value of $\bfV_{0,\widehat{\bfy}}^\star\bfT_{0,\bfy}$ is bounded uniformly from below for $\bbP$-almost all $y\in Y$.

\subsubsection{The special case of a single eigenvalue $m=1$}\label{sec:Eulerianrates}
The discrete Eulerian param\-eter-to-solution map $\bfscrE$ can be simplified for $m=1$. In particular, $\bfB_\bfy$ is a scalar and, since $\bfB_\bfy$ has the same eigenvalues as $\bfA_\bfy$ has eigenvalues in $\Lambda_c$, it holds that $\bfB_\bfy=\lambda_{0,\bfy}$ for $\bbP$-almost all $\bfy\in Y$. This simplifies the continuous and discrete Eulerian parameter-to-solution map to
\begin{align*}
\scrE&\colon Y\to\bbC\times X,
&&
\bfy\mapsto\scrE(\bfy)=(\lambda_{0,\bfy},u_{0,\widehat{\bfy}}),\\
\bfscrE^{(n)}&\colon Y\to\bbC\times\bbC^{N_n},
&&
\bfy\mapsto\bfscrE^{(n)}(\bfy)=(\lambda_{0,\bfy,n},\bfV_{0,\widehat{\bfy}})
\end{align*}
and implies that $\Lambda_c\cap\sigma(\bfA_\bfy)=\{\lambda_{0,\bfy}\}$ for $\bbP$-almost all $\bfy\in Y$. Thus, the convergence theory from \cref{sec:galerkinrates} applies and yields that
\begin{align*}
\|\scrE(\bfy)-\bfscrE^{(n)}(\bfy)\|_{\bbC\times X}
&\leq_{\phantom{\bfy,\widehat{\bfy}}}
|\lambda_{0,\bfy}-\lambda_{0,\bfy,n}|+\|u_{0,\widehat{\bfy}}-u_{0,\widehat{\bfy},n}\|_X\\
&\lesssim_{\bfy,\widehat{\bfy}}
|\lambda_{0,\bfy}-\lambda_{0,\bfy,n}|
+
|\lambda_{0,\widehat{\bfy}}-\lambda_{0,\widehat{\bfy},n}|\\
&\phantom{{}\lesssim|\lambda_{0,\bfy}-\lambda_{0,\bfy,n}|}+
\max_{\substack{x_0\in\ker(\calA_{\widehat{\bfy}}(\lambda_{0,\widehat{\bfy}}))\\\|x_0\|\leq 1}}\inf_{x_n\in X_n}\|x_0-x_n\|_X\\
&\lesssim_{\bfy,\widehat{\bfy}}
\delta_{n,\bfy}\delta_{n,\bfy}^\star+\delta_{n,\widehat{\bfy}}\delta_{n,\widehat{\bfy}}^\star+\max_{\substack{x_0\in\ker(\calA_{\widehat{\bfy}}(\lambda_{0,\widehat{\bfy}}))\\\|x_0\|\leq 1}}\inf_{x_n\in X_n}\|x_0-x_n\|_X,
\end{align*}
where the subscripts $\bfy$ and $\widehat{\bfy}$ indicate dependence of the implicit constants and $\delta_n$ and $\delta_n^\star$ on the parameter $\bfy$ and the reference value $\widehat{\bfy}$, respectively. Now, observe that the assumptions from \cref{sec:randomnonlinearevp} guarantee that the implicit constants in the above estimate are uniformly bounded for $\bbP$-almost all $\bfy\in Y$ and that $G(\calA_{\widehat{\bfy}},\Lambda_c)=\ker(\calA_{\widehat{\bfy}}(\lambda_{0,\widehat{\bfy}}))$. Introducing the error measure
\[
\eta_n
=
\esssup_{\bfy\in Y}\delta_{n,\bfy},
\qquad
\eta_n^\star
=
\esssup_{\bfy\in Y}\delta_{n,\bfy}^\star
\]
we obtain with \cref{eq:delta} that
\[
\|\scrE-\bfscrE^{(n)}\|_{L_\bbP^\infty(Y;\bbC\times X)}
\lesssim
\varepsilon_n
\coloneqq
\eta_n(\eta_n^\star+1).
\]
\begin{remark}
    The above estimates further reveal that the convergence rate improves to $\varepsilon_n=\eta_n\eta_n^\star$ when only the first component of the parameter-to-solution map is considered, i.e., when only the eigenvalue is of interest.
\end{remark}

We close our remarks by noting that the error analysis for $m>1$ is more involved mainly for two reasons. First, the interpretation of the eigenmatrix $\bfB_\bfy$ is more involved, and requires additional diagonalization before further steps. Second, also the convergence theory from \cref{sec:galerkinrates} is not necessarily applicable anymore, since the eigenvalues may bifurcate and the requirement $|\sigma(\calA_\bfy)\cap\Lambda_c|=1$ for $\bbP$-almost all $\bfy\in Y$ may not be satisfied anymore. Thus, for $m>1$, a convergence theory for clusters of eigenvalues and for the associated generalized eigenspace, would be required.
\newpage
\subsection{Lagrangian parameter-to-eigenspace maps}\label{sec:parameter-to-eigenspace-V2}
\subsubsection{Continuous perspective}
Our second approach for defining a parameter-to-solution map is based on a Lagrangian perspective. That is, we consider both, the eigenspace and the corresponding eigenmatrix, as functions of the parameter $\bfy\in Y$. We follow a similar strategy as it was done in \cite{GSH2023,HL2019b} for linear eigenvalue problems, and consider a basis of the generalized eigenspace $G(\calA_\bfy,\Lambda_c)$ obtained by applying $\calP_{G(\calA_\bfy,\Lambda_c)}$, cf.\ \cref{sec:Beynoriginal}. In \cite{GSH2023}, a thorough regularity analysis of the parameter-to-eigenspace map was provided for the specific map $\bfy\mapsto\calP_{G(\calA_\bfy,\Lambda_c)}\bfu_{0,\widehat{\bfy}}$, where $\bfu_{0,\widehat{\bfy}}$ is a basis of the generalized eigenspace $G(\calA_{\widehat{\bfy}},\Lambda_c)$ for a fixed reference value $\widehat{\bfy}\in Y$. We combine this perspective with Beyn's contour integral approach to define
\[
\bft_{0,\bfy}
\coloneqq
\calP_{G(\calA_\bfy,\Lambda_c)}\bfx,
\]
where $\bfx\in X^m$ is a vector with linearly independent entries in $X$ such that $\dim\operatorname{span}\bft_{0,\bfy}=m$ for $\bbP$-almost all $\bfy\in Y$. This implies that vector $\bft_{0,\bfy}$ is a basis of the generalized eigenspace $G(\calA_\bfy,\Lambda_c)$ for $\bbP$-almost all $\bfy\in Y$. In accordance with Beyn's approach, the vector $\bfx$ may be realized by drawing $m$ random samples from $X$ according to a probability distribution on $X$. Now, let $\bfu_\bfy$ be an orthonormal basis of generalized eigenelements of $\operatorname{span}(\bft_{0,\bfy})=G(\calA_\bfy,\Lambda_c)$ corresponding to the eigenvalues in $\Lambda_c$ and $\bfLambda_\bfy$ a Jordan matrix such that
\[
\calA_\bfy(\bfLambda_\bfy)\bfu_\bfy=\bfNull.
\]
Then, there is a nonsingular matrix $\bfR_\bfy$ such that $\bft_{0,\bfy}\bfR_\bfy=\bfu_\bfy$. Setting $\bfC_\bfy=\bfR_\bfy\bfLambda_\bfy\bfR_\bfy^{-1}$ yields
\[
\calA_\bfy(\bfC_\bfy)\bft_{0,\bfy}=\bfNull.
\]
Thus, we define our continuous Lagrangian parameter-to-solution map as
\[
\scrL\colon Y\to\bbC^{m\times m}\times X^m,
\qquad
\bfy\mapsto\scrL(\bfy)\coloneqq(\bfC_\bfy,\bft_{0,\bfy}).
\]
The assumptions from \cref{sec:randomnonlinearevp} guarantee that
\[
\scrL
\in
L_\bbP^\infty(Y;\bbC^{m\times m}\times X^m).
\]
whenever the $m$-th singular value of $\bfT_{0,\bfy}$ is bounded uniformly from below for $\bbP$-almost all $y\in Y$.

\subsubsection{Discretization and implementation}
We consider again a Galerkin discretization of \cref{eq:randomnonlinearevp} as in \cref{sec:galerkin} and the integrals \cref{eq:T01y} from Beyn's contour integral approach, but follow the derivation in \cref{sec:Beynmodified} to define
\[
\bfC_\bfy
=
\bfT_{0,\bfy}^\dagger\bfT_{1,\bfy}.
\]
Again, the nonzero eigenvalues of $\bfC_\bfy$ coincide with the eigenvalues of $\bfA_\bfy$ in $\Lambda_c$, including algebraic multiplicities and with at least $m'-m$ eigenvalues of $\bfC_\bfy$ vanishing.
Our discrete parameter-to-eigenspace map is then given as
\[
\bfscrL\colon Y\to\bbC^{m'\times m'}\times\bbC^{N_n\times m'},
\qquad
\bfy\mapsto\bfscrL(\bfy)
\coloneqq
(\bfC_\bfy,\bfT_{0,\bfy}).
\]
Understood in the sense of \cref{sec:NEVPeigenspaces}, $\bfscrL(\bfy)$ is a solution to the eigenvalue problem
\[
\bfA_\bfy(\bfC_\bfy)\bfT_{0,\bfy}=0,
\]
for $\bbP$-almost every $\bfy\in Y$. In contrast to the Eulerian parameter-to-eigenspace map $\bfscrE$, the Lagrangian map $\bfscrL(\bfy)$ provides a direct insight into the behaviour of the eigenspaces. The discrete parameter-to-eigenspace map $\bfscrL$ is well-defined for $\bbP$-almost all $\bfy\in Y$ and it holds
\[
\bfscrL\in L_\bbP^\infty\Big(Y;\bbC^{m'\times m'}\times\bbC^{N_n\times m'}\Big).
\]

\subsubsection{The special case of a single eigenvalue $m=1$}\label{sec:Lagrangianrates}
By the same argumentation as for the Eulerian parameter-to-eigenspace map, the discrete Lagrangian parameter-to-eigenspace map $\bfscrL$ can be simplified for $m=m'=1$. In particular, $\bfC_\bfy$ is a scalar and, since $\bfC_\bfy$ has the same eigenvalues as $\bfA_\bfy$ has eigenvalues in $\Lambda_c$, it holds that $\bfC_\bfy=\lambda_{0,\bfy}$ for $\bbP$-almost all $\bfy\in Y$. This simplifies the discrete Lagrangian parameter-to-eigenspace map to
\begin{align*}
\scrL&\colon Y\to\bbC\times X,
&&
\bfy\mapsto\scrL(\bfy)=(\lambda_{0,\bfy},t_{0,\bfy}),\\
\bfscrL^{(n)}&\colon Y\to\bbC\times\bbC^{N_n},
&&
\bfy\mapsto\bfscrL^{(n)}(\bfy)=(\lambda_{0,\bfy,n},\bfT_{0,\bfy}).
\end{align*}
In analogy to the Eulerian case, convergence theory from \cref{sec:galerkinrates} and the assumptions from \cref{sec:randomnonlinearevp} yield that
\begin{align*}
\|\scrL(\bfy)-\bfscrL^{(n)}(\bfy)\|_{\bbC\times X}
&\leq_{\phantom{\bfy,\widehat{\bfy}}}
|\lambda_{0,\bfy}-\lambda_{0,\bfy,n}|+\|t_{0,\bfy}-t_{0,\bfy,n}\|_X\\
&\lesssim_{\bfy,\widehat{\bfy}}
|\lambda_{0,\bfy}-\lambda_{0,\bfy,n}|
+
\max_{\substack{x_0\in\ker(\calA_\bfy(\lambda_{0,\bfy}))\\\|x_0\|\leq 1}}\inf_{x_n\in X_n}\|x_0-x_n\|_X\\
&\lesssim_{\bfy,\widehat{\bfy}}
\delta_{n,\bfy}\delta_{n,\bfy}^\star+\max_{\substack{x_0\in\ker(\calA_\bfy(\lambda_{0,\bfy}))\\\|x_0\|\leq 1}}\inf_{x_n\in X_n}\|x_0-x_n\|_X\\
&\leq_{\phantom{\bfy,\widehat{\bfy}}}
\eta_n(\eta_n^\star+1)
=
\varepsilon_n.
\end{align*}
This yields
\[
\|\scrL-\bfscrL^{(n)}\|_{L_\bbP^\infty(Y;\bbC\times X)}
\lesssim
\varepsilon_n.
\]
In analogy to the Eulerian case, our estimates show that the convergence rate improves to $\varepsilon_n=\eta_n\eta_n^\star$ when only the eigenvalue is of interest. Moreover, the same remarks as for the Eulerian parameter-to-eigenspace map apply for the Lagrangian parameter-to-eigenspace map for $m>1$.

\subsection{Multilevel hierarchies of parameter-to-solution maps}\label{sec:multilevelhierarchies}
One of the challenges for multilevel algorithms for eigen\-spaces and parametrized eigenvalue problems, even in the linear case, is to consistently relate eigenspaces across different discretization levels and parameter values. More precisely, one needs to relate a basis or eigenmatrix of the generalized eigenspace $G(\bfA_\bfy,\Lambda_c)$ for $\bbP$-almost every $\bfy\in Y$ both across different choices of $\bfy\in Y$ and across different discretization levels.

For a fixed discretization, our discrete parameter-to-eigenspace maps $\bfscrE$ and $\bfscrL$ elegantly circumvent this issue by anchoring the basis to the previously chosen matrix $\widehat{\bfV}$, which remains fixed for all $\bfy\in Y$. Once $\widehat{\bfV}$ is fixed, $\bfscrE$ and $\bfscrL$ are uniquely determined for $\bbP$-almost all $\bfy\in Y$, completely removing the need to match eigenspaces across different parameter values.

To extend this approach to multilevel hierarchies of parameter-to-eigenspace maps, we consider discretizations of the nonlinear eigenvalue problem as in \cref{sec:galerkin} and fix $m'\geq m$ random elements $\widehat{\bfv}=\{v_1,\ldots,v_{m'}\}\subset X$. Following the philosophy of \cite{Bey2012}, we may assume that $\calP_{G(\calA_\bfy,\Lambda_c)}\widehat{\bfv}\neq\bfNull$. With $\widehat{\bfv}$ fixed, we define
\begin{equation}\label{eq:Vhat-disc}
\widehat{\bfV}^{(n)}
=
\bigg[
[\langle v_1,\phi_i^{(n)}\rangle_X]_{i=1}^{N_n},\ldots,[\langle v_{m'},\phi_i^{(n)}\rangle_X]_{i=1}^{N_n}
\bigg]
\in\bbC^{N_n\times m'},
\end{equation}
implying that $\bfA_{\bfy}^{(n)}(z)^{-1}\widehat{\bfV}^{(n)}$ is well-defined for all $z\in\Gamma$ and $\bbP$-almost all $\bfy\in Y$. In particular, $\bfA_{\bfy}^{(n)}(z)^{-1}\widehat{\bfV}^{(n)}$ is a matrix whose columns are the coefficients of Galerkin approximations of the solutions to the linear problems $\calA_\bfy(z)^{-1}v_j$ for $j=1,\ldots,m'$. Well-definedness of these problems and Galerkin approximations follows from the assumptions in \cref{sec:randomnonlinearevp}. In particular, our assumptions guarantee that the Galerkin solutions depend continuously on $\widehat{\bfv}$ and that the dependence is uniform for $\bbP$-almost all $\bfy\in Y$ and for all $z\in\Gamma$.

In complete analogy to the case of a fixed discretization, introduce
\begin{align*}
		\bfT_{0,\bfy}^{(n)}
		&\coloneqq\calP_{G(\bfA_\bfy^{(n)},\Lambda_c)}\widehat{\bfV}^{(n)}\coloneqq
		\frac{1}{2\pi i}\int_{\Gamma}\bfA_\bfy^{(n)}(z)^{-1}\widehat{\bfV}^{(n)}\dd z\in\bbC^{N_n\times m'},\\
		\bfT_{1,\bfy}^{(n)}
		&\coloneqq\phantom{\coloneqq\calP_{G(\bfA_\bfy^{(n)},\Lambda_c)}\widehat{\bfV}^{(n)}}
		\frac{1}{2\pi i}\int_{\Gamma}z\bfA_\bfy^{(n)}(z)^{-1}\widehat{\bfV}^{(n)}\dd z\in\bbC^{N_n\times m'},
\end{align*}
and the discrete parameter-to-eigenspace maps
\[
\bfy\mapsto\bfscrE^{(n)}(\bfy)\coloneqq(\bfB_\bfy^{(n)},\bfV_0^{(n)}),
\qquad
\bfy\mapsto\bfscrL^{(n)}(\bfy)\coloneqq(\bfC_\bfy^{(n)},\bfT_{0,\bfy}^{(n)}),
\]
where $\bfB_\bfy^{(n)}$ and $\bfC_\bfy^{(n)}$ are defined in complete analogy to \cref{sec:parameter-to-eigenspace-V1,sec:parameter-to-eigenspace-V2} and $\bfV_0^{(n)}$ is a reference basis of $G(\bfA_{\widehat{\bfy}}^{(n)},\Lambda_c)$ for a fixed reference value $\widehat{\bfy}\in Y$. For $m>1$, the left singular vectors of $\bfT_{0,\widehat{\bfy}}^{(n)}$ are only unique up to a unitary matrix $\bfU_n$ and choose $m'=m$ and $\bfV_0^{(n)}=\bfT_{0,\widehat{\bfy}}^{(n)}$. This implies that the discrete parameter-to-eigenspace maps $\bfscrE^{(n)}$ and $\bfscrL^{(n)}$ are well-defined for $\bbP$-almost all $\bfy\in Y$ and for all $n\in\bbN$. In particular, $\bfscrE^{(n)}(\bfy)$ and $\bfscrL^{(n)}(\bfy)$ are uniquely determined for $\bbP$-almost all $\bfy\in Y$ once $\widehat{\bfv}$ is fixed. This allows us to define multilevel hierarchies of parameter-to-eigenspace maps which are well-defined for $\bbP$-almost all $\bfy\in Y$.

\section{Multilevel Monte Carlo methods}\label{sec:MLMC}

\subsection{Plain vanilla Monte Carlo}
Having well-defined parameter-to-eigenspace maps $\scrE$ and $\scrL$ as well as their discretizations at our disposal, we are now in the position to compute statistical quantities of interest therefrom. We consider the means
\[
\bbE[\scrS]
=
\int_{\bfy}\scrS(\bfy)\dd\bbP(\bfy),
\qquad
\scrS\in\{\scrE,\scrL\},
\]
of the Eulerian and Lagrangian parameter-to-solution maps as exemplary and representative target statistical quantities of interest. For the computational approximation of the integral we consider the Monte Carlo method combined with a Galerkin discretization of the parameter-to-solution maps, i.e.,
\[
\bbE[\scrS]
\approx
\frac{1}{M}\sum_{j=1}^M\bfscrS^{(n)}(\bfy_j)
\eqqcolon
E_M[\bfscrS^{(n)}],
\qquad
\scrS\in\{\scrE,\scrL\},
\]
where $\bfy_j$ are iid random variables drawn from the probability measure $\bbP$ and $\bfscrS^{(n)}$ is the discrete approximation of $\scrS$ using $X_n\subset X$. It is well known that the Monte Carlo method converges with rate $M^{-1/2}$ in the root mean square sense. Moreover, for $m=1$ and in accordance with common Monte Carlo approximation theory, cf.\ e.g., \cite{BSZ2011}, and \cref{sec:Eulerianrates,sec:Lagrangianrates}, we obtain the error estimates
\begin{align*}
&\big\|\bbE[\scrS]-E_M\big[\bfscrS^{(n)}\big]\big\|_{L_\bbP^2(Y;\bbC\times X)}\\
&\qquad\leq
\big\|\bbE[\scrS]-\bbE\big[\bfscrS^{(n)}\big]\big\|_{L_\bbP^2(Y;\bbC\times X)}
+
\big\|\bbE[\bfscrS^{(n)}]-E_M\big[\bfscrS^{(n)}\big]\big\|_{L_\bbP^2(Y;\bbC\times X)}\\
&\qquad\lesssim
\varepsilon_n+M^{-1/2}\|\bfscrS^{(n)}\|_{L_\bbP^2(Y;\bbC\times X)}\\
&\qquad\leq
\varepsilon_n+M^{-1/2}(\|\scrS\|_{L_\bbP^2(Y;\bbC\times X)}+\varepsilon_n)\\
&\qquad\lesssim
\varepsilon_n+M^{-1/2},
\qquad
\scrS\in\{\scrE,\scrL\}.
\end{align*}
It is clear that $M^{-1/2}\sim\varepsilon_n$ balances the error between Monte Carlo and Galerkin approximation. For most cases, reaching this balance requires a large number of Monte Carlo samples $M$ and, thus, a large number of discrete nonlinear eigenvalue solves with computational cost $N$. The resulting overall cost of $M\cdot N$ is often computationally prohibitive for $M$ and $N$ large.

\subsection{Multilevel Monte Carlo}
The multilevel Monte Carlo method, cf., e.g., \cite{Gil2015,Hei2001}, is a variance reduction technique which allows to reduce the number of samples on the finest discretization level. The idea is to consider a series of discretizations $X_0,X_1,\ldots,X_L$ of increasing accuracy and to write the mean of the parameter-to-eigenspace maps on the finest discretization level $X_L$ as a telescopic sum
\[
\begin{aligned}
\bbE[\bfscrS^{(L)}]
&=
\bbE[\bfscrS^{(0)}]
+
\sum_{\ell=1}^L\bbE[\bfscrS^{(\ell)}-\bfscrS^{(\ell-1)}]\\
&\approx
E_{M_0}[\bfscrS^{(0)}]
+
\sum_{\ell=1}^L E_{M_\ell}[\bfscrS^{(\ell)}-\bfscrS^{(\ell-1)}]
\eqqcolon
E^{ML}[\bfscrS^{(L)}],
\qquad
\bfscrS\in\{\bfscrE,\bfscrL\},
\end{aligned}
\]
For $m=1$ and in accordance with common multilevel Monte Carlo approximation theory, cf.\ e.g., \cite{BSZ2011,Gil2015}, and \cref{sec:galerkinrates}, we obtain the error estimate
\begin{equation*}
    \begin{split}
    \big\|\bbE[\scrS]-E^{ML}[\bfscrS^{(L)}]\big\|_{L_\bbP^2(Y;\bbC\times X)}
    \lesssim
    \varepsilon_L
    &+
    M_0^{-1/2}\|\bfscrS^{(0)}\|_{L_\bbP^2(Y;\bbC\times X)} \\
    &+
    \sum_{\ell=1}^L M_\ell^{-1/2}(\varepsilon_\ell+\varepsilon_{\ell-1}).
    \end{split}
\end{equation*}
Here, we used that $\|\bfscrS^{(\ell)}-\bfscrS^{(\ell-1)}\|_{L_\bbP^2(Y;\bbC\times X)}\lesssim\varepsilon_\ell+\varepsilon_{\ell-1}$ for $\ell\geq 1$. Assuming that $\varepsilon_\ell\lesssim\varepsilon_{\ell-1}$, choosing $M_0\sim\varepsilon_L^{-2}$ and $M_\ell\sim(\varepsilon_{\ell-1}/\varepsilon_L)^2$ for $\ell\geq 1$ balances the error between Monte Carlo and Galerkin approximation such that the overall error is of order $(L+1)\varepsilon_L$. This shifts most of the required samples to coarser discretization levels, with significantly lower cost than on the finest discretization level. The resulting overall cost of the multilevel Monte Carlo method is then given by $\sum_{\ell=0}^L M_\ell N_\ell$. The achievable cost savings depend on the convergence rates of the Galerkin discretization and the cost of a single nonlinear eigenvalue solve on each discretization level, cf.\ e.g., \cite{Gil2015,Hei2001}.

\begin{remark}
The approach is straightforwardly adapted to other quadrature methods such as quasi-Monte Carlo and sparse grids. For quadratures of order $M_\ell^{-\alpha}$, $\alpha>0$, this leads to sample numbers $M_\ell\sim(\varepsilon_{\ell}/\varepsilon_L)^{1/\alpha}$. A rigorous analysis requires the parametric regularity which is out of the scope of the present article.
\end{remark}
\begin{remark}
Most applications of the multilevel Monte Carlo methods require a hierarchy of discretization spaces $X_0\subset X_1\subset\ldots\subset X_L$ and corresponding prolongation operators $P_\ell\colon X_\ell\to X_{\ell+1}$ for $\ell=0,\ldots,L-1$. For the setting of this article, this is only required if one is interested in the eigenspaces of the Lagrangian parameter-to-solution map. For the Eulerian parameter-to-solution map and the eigenmatrix of the Lagrangian map, the multilevel Monte Carlo method can be applied without the need of a hierarchy of discretization spaces and specifically without the need of prolongation operators, as the algorithm only needs to operate on the eigenmatrices $\bfB_\bfy^{(\ell)}$ and $\bfC_\bfy^{(\ell)}$ on each discretization level $\ell=0,\ldots,L$. The size of those matrices is independent of the discretization level, avoiding the need for prolongation operators. This requires the same $\widehat{\bfv}$ on all levels and, for the Eulerian map, a level-consistent $\bfV_0^{(\ell)}$, see \cref{sec:multilevelhierarchies}.
\end{remark}
\begin{remark}\label{rem:qoi}
For $m>1$, the eigenvalues of $\bbE[\bfB_\bfy]$ or $\bbE[\bfC_\bfy]$ are in general not the means of the eigenvalues of $\calA_\bfy$ in $\Lambda_c$, and need not be real. Similarity invariant quantities such as $\operatorname{tr}(\bfB_\bfy^k)=\operatorname{tr}(\bfC_\bfy^k)=\sum_{i=0}^{m-1}\lambda_{i,\bfy}^k$, $k\geq 1$, are unaffected by crossings and can be estimated verbatim.
\end{remark}

\begin{algorithm}
  \DontPrintSemicolon
  \SetKwFunction{FMC}{PlainMC}
  \SetKwFunction{FMLMC}{MLMC}
  \SetKwProg{Fn}{Function}{:}{\KwRet}
  \Fn{\FMC{$\bfy \mapsto (\bfA_{\bfy}(\bullet) \colon \Lambda \to \bbC^{N \times N} )$, $\widehat{\bfV} \in \bbC^{N \times m'}$, $\{\bfy_j\}_{j=1}^M$}}{
        \For{$j=1,\ldots, M$}{
            Compute $\bfT_{0, \bfy_j}, \bfT_{1, \bfy_j} \in\bbC^{N\times m'}$ as in \eqref{eq:T01y}\;
            Compute $\bfC_{\bfy_j} = \bfT_{0, \bfy_j}^\dagger \bfT_{1, \bfy_j} \in \bbC^{m' \times m'}$\;
        }
        \Return{$\overline{\bfC} = \frac{1}{M} \sum_{j=1}^M \bfC_{\bfy_j}, \overline{\bfT}_{0} = \frac{1}{M} \sum_{j=1}^M \bfT_{0, \bfy_j}$}
  }
  \SkipBeforeAndAfter
  \SetKwProg{Pn}{Function}{:}{\KwRet}
  \Pn{\FMLMC{{$\{\bfy \mapsto (\bfA_\bfy^{(\ell)}(\bullet) \colon \Lambda \to \bbC^{N_\ell \times N_\ell})\}_{\ell = 0}^L$, $\{M_\ell\}_{\ell = 0}^L$}}}{
        Generate $m' \geq m$ random elements $\widehat{\bfv} = \{v_1, \ldots, v_{m'}\}$\;
        \For{$\ell = 0,\ldots,L$}{
            Compute the discretization $\widehat{\bfV}^{(\ell)} \in \bbC^{N_\ell \times m'}$ as in \eqref{eq:Vhat-disc}\;
            Generate samples $\bfy_j \sim \bbP, \: j=1,\ldots, M_\ell$\;
            \makebox[0pt][l]{$\bfC^{(\ell), M_\ell},$}\phantom{$\bfC^{(\ell-1), M_\ell},$} \makebox[0pt][l]{$\bfT_0^{(\ell), M_\ell}$}\phantom{$\bfT_0^{(\ell-1), M_\ell}$} = \FMC{\makebox[0pt][l]{$\bfA^{(\ell)}_\bfy$,}\phantom{$\bfA^{(\ell-1)}_\bfy$,} \makebox[0pt][l]{$\widehat{\bfV}^{(\ell)}$,}\phantom{$\widehat{\bfV}^{(\ell-1)}$,} $\{\bfy_j\}_{j=1}^{M_\ell}$}\;
            $\bfC^{(\ell-1), M_\ell}$, $\bfT_0^{(\ell-1), M_\ell}$ = \FMC{$\bfA^{(\ell-1)}_\bfy$, $\widehat{\bfV}^{(\ell-1)}$, $\{\bfy_j\}_{j=1}^{M_\ell}$}\;
            Set $\bfD^{(\ell)} = \bfC^{(\ell), M_\ell} - \bfC^{(\ell-1), M_\ell}$\;
            Set $\boldsymbol{\Delta}^{(\ell)} = \bfT^{(\ell), M_\ell}_{0} - \bfI_{\ell-1}^\ell\bfT^{(\ell-1), M_\ell}_{0}$\tcp*{$\bfI_{\ell-1}^\ell$ prolongation matrix}
        }
        \Return{$\bfC^{\text{ML}} = \sum_{\ell=0}^L \bfD^{(\ell)}$, $\bfT_0^{\text{ML}} = \sum_{\ell=0}^L \bfI_\ell^L \boldsymbol{\Delta}^{(\ell)}$}    
  }
  \caption{Multilevel Monte-Carlo (MLMC) quadrature for the Lagrangian parameter-to-eigenspace map}\label{algo:mlmc}
\end{algorithm}

\section{Numerical experiments: random diffusion equation}\label{sec:experiments_toy}
\subsection{The random diffusion eigenvalue problem}
We consider the eigenvalue problem for the diffusion operator with random coefficient and with homogeneous boundary conditions, i.e.,
\begin{subequations}\label{eq:DiffusionEVP}
\begin{align}
-\divergence (a_\bfy(\bullet) \nabla u_\bfy)&=\lambda_\bfy u_\bfy,&&\text{in}~D\coloneqq [0,1]^2,\\
\qquad
u_\bfy&=0,&&\text{on}~\partial D,
\end{align}
\end{subequations}
where the diffusion coefficient $a_\bfy$ is given by a Karhunen-Lo\`eve expansion
\[
	a_\bfy (\bfx) = \bar{a}(\bfx) + \sum_{j=1}^{20} \psi_j(\bfx) y_j\in C(\overline{D}),
\]
with $\bar{a}(\bfx) \equiv 1$, $\bfy = [y_1, \ldots, y_{20}]^T \in [-\frac{1}{3},\frac{1}{3}]^{20}$, and $\psi_j(\bfx) = 2^{-j} \sin(j\pi x_1)\sin(j\pi x_2)$. Clearly, the diffusion coefficient is bounded uniformly from below and above, implying that \cref{eq:DiffusionEVP} is a linear, elliptic eigenvalue problem. Its nonlinear formulation compatible with \cref{eq:nonlinearevp} is
\[
\calA(\lambda_\bfy)u_\bfy
\coloneqq
-\divergence (a_\bfy(\cdot) \nabla u)-\lambda u
=
0.
\]

It is well known that, for $\bfy = 0$, the exact eigenpairs of \cref{eq:DiffusionEVP} are given by
\[
    \lambda_\bfk = \pi^2(k_1^2+k_2^2)
\]
with associated eigenfunctions
\[
    u_\bfk(\bfx) = \sin(\pi k_1 x_1)\sin(\pi k_2 x_2)
\]
for all $\bfk =(k_1,k_2) \in \bbN^2$.

\subsection{Computational setup}
For the Galerkin discretization of \cref{eq:DiffusionEVP} we use nested, continuous, piecewise linear finite element spaces arising from uniform refinement of an initial triangular mesh. The discretized version of \cref{eq:DiffusionEVP} thus reads
\begin{equation}\label{eq:DiffusionEVP-discr}
	\bfA^{(n)}_\bfy(\lambda^{(n)}_\bfy)\bfu_\bfy \coloneqq (\bfS^{(n)}_\bfy - \lambda^{(n)}_\bfy \bfM^{(n)}) \bfu_\bfy = 0,
\end{equation}
where 
\begin{align*}
	\bfS^{(n)}_\bfy &= \left[(a_\bfy(\bullet)\nabla \phi_j^{(n)}, \nabla \phi_i^{(n)})_{L^2(D)}\right]_{i,j=1,\ldots,N^{(n)}},\\
	\bfM^{(n)} &= \left[(\phi_j^{(n)}, \phi_i^{(n)})_{L^2(D)}\right]_{i,j=1,\ldots,N^{(n)}},
\end{align*}
are the stiffness and mass matrix, respectively, and $\bfu_\bfy$ is the coefficient vector of a function in $V^{(n)} = \operatorname{span}\{\phi_1^{(n)}, \ldots, \phi_{N^{(n)}}^{(n)}\} \subset H^1_0(D)$.

We choose a circle-shaped contour with center $\mu = 38$ and radius $r = 28$, which contains the two smallest distinct eigenvalues $2\pi^2$ and $5\pi^2$ in the reference case $\bfy = 0$, i.e., $m=3$ eigenvalues counting algebraic multiplicities. We observe that while in the reference case the second eigenvalue is double, it bifurcates for other values of $\bfy$, resulting in two distinct simple eigenvalues. 

As the finite element method will give us a convergence rate of $\mathcal{O}(h^2)$ for the approximation of the eigenvalues, where $h$ is the mesh size, in our multilevel quasi-Monte Carlo method the number of samples has to be adapted according to $\mathcal{O}(h^{-2})$.
For the quasi-Monte Carlo quadrature we use the Kronecker and the Halton sequence \cite{Dick_Kuo_Sloan_2013}. The number of samples per level are given in \cref{tab:number-of-samples_toy}.

\begin{table}[htb]
    \centering
\begin{tabular}{c c c c c c c c}
    \hline
    $\ell=2$ & $\ell=3$ & $\ell=4$ & $\ell = 5$ & $\ell=6$ & $\ell=7$ & $\ell = 8$ & $\ell = 9$ \\\hline
    163\,840 & 40\,960  & 10\,240  & 2\,560     & 640      & 160      & 40         & 10 \\\hline
    \end{tabular}
    \medskip
    \caption{Number of QMC samples per level for the diffusion problem with random diffusion coefficient.}
    \label{tab:number-of-samples_toy}
\end{table}
The computations were carried out on four nodes of a cluster, with each node equipped with two AMD EPYC 9354 32-Core processors each, which amounts to 256 CPU cores in total.

\subsection{Numerical results}\label{subsec:toyproblem-numerical-results}
To measure the approximation error of the parameter-to-eigenmatrix maps, we compute the eigenvalues $\lambda^L_1, \ldots, \lambda^L_{m}$ of the multilevel approximations $E^{ML}[\bfB^{(L)}]$ and $E^{ML}[\bfC^{(L)}]$ and compute their error to the reference values $\lambda^{\text{ref}}_1, \ldots, \lambda^{\text{ref}}_{m}$ in the $\ell^\infty$-distance. That is, we compute $\max_{k=1,\ldots,m} |\lambda^L_k - \lambda^{\text{ref}}_k|$. We emphasize that, for $m>1$, the eigenvalues of $\bbE[\bfB_\bfy]$ and $\bbE[\bfC_\bfy]$ are, in general, not the expectations of the individual eigenvalues of $\calA_\bfy$ in $\Lambda_c$.

To measure the error for the Lagrangian parameter-to-eigenspace maps we also require an error measure for the matrices $E^{ML}[\bfT_{0,\bfy}^{(L)}]$, i.e., the approximations to $\bbE[\calP_{G(\calA_\bfy,\Lambda_c)}]\widehat{\bfv}$, which we refer to as expected eigenspaces in the following.
To this end, we first orthonormalize the columns of $E^{ML}[\bfT_{0,\bfy}^{(L)}]$, computed by means of a different quadrature method, with respect to the $L^2(D)$-inner product to obtain a reference $L^2(D)$-orthonormal basis $\bfT^\text{ref}$. We then compute
\[
\Big\lVert (\Id-\calP_{L^2(D), \bfT^{\text{ref}}})E^{ML}[\bfT_{0,\bfy}^{(L)}]\Big\rVert_{[L^2(D)]^m},
\]
where $\calP_{L^2(D), \bfT^{\text{ref}}}$ denotes the projection with respect to the $L^2(D)$-inner product onto the subspace spanned by the columns of $\bfT^{\text{ref}}$.

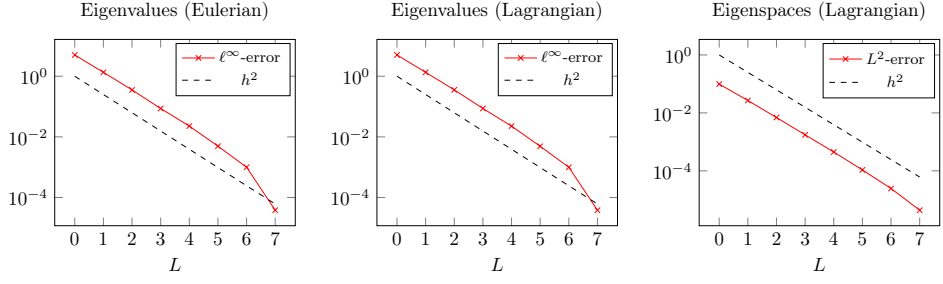
\begin{figure}
    \centering
    \begin{tikzpicture}[scale=0.72]
    \begin{groupplot}[
        group style={
            group size=3 by 1,
            horizontal sep=1.5cm,
        },
        width=6cm,
        height=5cm,
        ymode=log,
        xlabel={$L$},
        xtick={0,1,2,3,4,5,6,7},
        legend style={font=\footnotesize}
    ]
    \nextgroupplot[title={Eigenvalues (Eulerian)}]
    \addplot[red, mark=x]table[x=L,y=eigenvalues_euler]{figures/diff_coeff_error.data}; \addlegendentry{$\ell^\infty$-error}
    \addplot[domain=0:7, samples=200, dashed]{2^(-2*x)}; \addlegendentry{$h^2$}

    \nextgroupplot[title={Eigenvalues (Lagrangian)}]
    \addplot[color=red,mark=x]table[x=L,y=eigenvalues_lagrange]{figures/diff_coeff_error.data};\addlegendentry{$\ell^\infty$-error}
    \addplot[domain=0:7, samples=200, dashed]{2^(-2*x)}; \addlegendentry{$h^2$}

    \nextgroupplot[title={Eigenspaces (Lagrangian)}]
    \addplot[color=red,mark=x]table[x=L,y=eigenspaces_lagrange]{figures/diff_coeff_error.data};\addlegendentry{$L^2$-error}
    \addplot[domain=0:7, samples=200, dashed]{2^(-2*x)}; \addlegendentry{$h^2$}
    
    \end{groupplot}
    \end{tikzpicture}

    \caption{Convergence of multilevel QMC for the Eulerian and Lagrangian parameter-to-solution maps for the diffusion problem with random diffusion coefficient.}
    \label{fig:plots-toy-prob}
\end{figure}

\begin{figure}
    \centering
    \includegraphics[width=\textwidth]{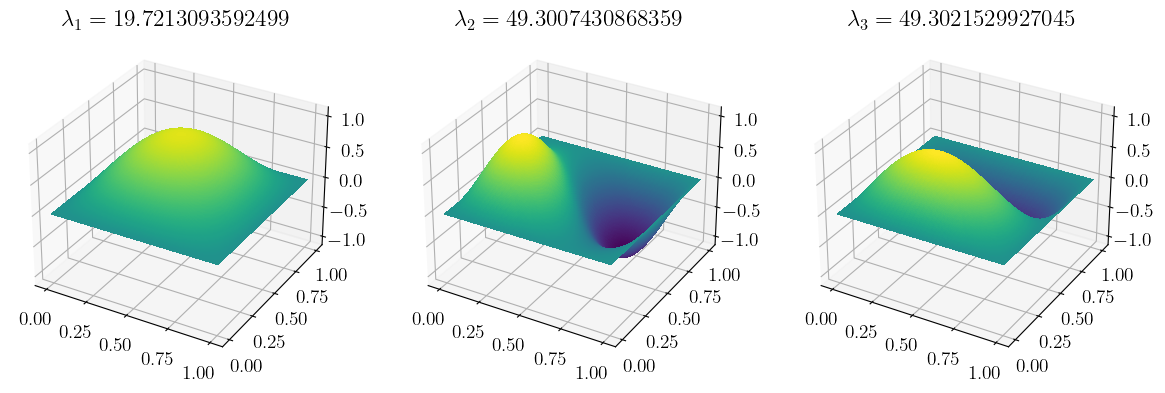}
    \caption{Visualization of the expected eigenfunctions (Lagrangian para\-meter-to-eigenspace map).}
    \label{fig:visual-2d-exp-eigfct}
\end{figure}

\Cref{fig:plots-toy-prob} shows the convergence of the multilevel methods, using the Kronecker sequence for the quasi-Monte Carlo quadrature, for both the Eulerian and Lagrangian parameter-to-solution maps. As a single-level Monte Carlo approach was out of reach, we used a multilevel quasi-Monte Carlo method with the Halton sequence to obtain a reference solution.
We observe that the convergence rates agree with the theoretical predictions for $m=1$ and we obtain a convergence rate of $h^2$ for the eigenspaces due to measuring the error in the $L^2$-norm as opposed to the $H^1$-norm, which would give us a rate of $h$. 
\Cref{fig:visual-2d-exp-eigfct} shows a visualization of the expected eigenfunctions extracted using the Lagrangian parameter-to-eigenspace map.

\section{Numerical experiments: Nonlinear formulation of the Laplacian eigenvalue problem}\label{sec:experiments_BEM}
\subsection{Boundary integral equations for the Laplacian eigenvalue problem}
We consider the eigenvalue problem of the Laplace operator with homogeneous boundary conditions, i.e.,
\begin{subequations}\label{eq:LaplaceEVP}
\begin{align}
-\Delta u_\bfy&=\lambda_\bfy u_\bfy,&&\text{in}~D_\bfy,\\
\qquad
u_\bfy&=0,&&\text{on}~\partial D_\bfy,
\end{align}
\end{subequations}
where $D_\bfy$ is a randomly parametrized, bounded, and simply connected Lipschitz domain $D_\bfy\subset\bbR^3$. While this eigenvalue problem can be transformed onto a fixed reference domain to obtain a linear eigenvalue problem with random coefficients as in \cref{eq:DiffusionEVP}, classical approaches for its numerical solution by Galerkin discretization require a volume mesh of the reference or for each $D_\bfy$ for $\bbP$-almost all $\bfy\in Y$. In the following, we circumvent the need for a volume mesh by following the approach of Steinbach and Unger \cite{SU2009,SU2012} to reformulate \cref{eq:LaplaceEVP} as a nonlinear eigenvalue problem on the boundaries $\{\partial D_\bfy\}_{\bfy\in Y}$. To achieve this, we set $\kappa_\bfy^2=\lambda_\bfy$ and reformulate \cref{eq:LaplaceEVP} in terms of the Helmholtz equation
\begin{subequations}\label{eq:Helmholtz}
\begin{align}
-\Delta u_\bfy-\kappa_\bfy^2 u_\bfy&=0,&&\text{in}~D_\bfy,\\
\qquad
u_\bfy&=0,&&\text{on}~\partial D_\bfy.
\end{align}
\end{subequations}
Then, any solution $u_\bfy\in H_0^1(D_\bfy)$ to \cref{eq:Helmholtz} can be represented as a potential
\begin{equation}\label{eq:potential}
u_\bfy(\bfx)
=
(\tilde{\calV}_\bfy(\kappa_\bfy)t_\bfy)(\bfx)
=
\frac{1}{4\pi}\int_{\partial D_\bfy}\frac{e^{i\kappa_\bfy\|\bfx-\bfz\|_2}}{\|\bfx-\bfz\|_2}t_\bfy(\bfz)\dd\sigma_{\bfz},
\qquad
\bfx\in D_\bfy,
\end{equation}
for some unknown $t_\bfy\in H^{-1/2}(\partial D_\bfy)$, cf.\ e.g., \cite{Ste2008}. Taking the interior trace $\gamma_0^{\text{int}}$ to \cref{eq:potential} yields
\[
0
=
(\gamma_0^{\text{int}}u_\bfy)(\bfx)
=
\big(\gamma_0^{\text{int}}\tilde{\calV}_\bfy(\kappa_\bfy)t_\bfy\big)(\bfx)
\defis
(\calV_\bfy(\kappa_\bfy)t_\bfy)(\bfx).
\]
It was shown in \cite{SU2009,SU2012} that if $(\kappa_\bfy^2,u_\bfy)\in\bbR\times H_0^1(D_\bfy)\setminus\{0\}$ is an eigenpair of \cref{eq:LaplaceEVP}, then $t_\bfy=\frac{\partial u_\bfy}{\partial\bfn_\bfy}$ on $\partial D_\bfy$ is a nontrivial solution to
\begin{equation}\label{eq:sloevp}
\calV_\bfy(\kappa_\bfy)(t_\bfy)=0
\qquad
\text{on}~\partial D_\bfy.
\end{equation}
Conversely, any solution $(\kappa_\bfy,t_\bfy)\in\bbR\times H^{-1/2}(\partial D_\bfy)\setminus\{0\}$ to \cref{eq:sloevp} yields an eigenpair $\big(\kappa_\bfy^2,\tilde{\calV}_\bfy(\kappa_\bfy)(t_\bfy)\big)\in\bbR\times H_0^1(D_\bfy)\setminus\{0\}$ of \cref{eq:LaplaceEVP}.

Thus, finding $(\kappa_\bfy,t_\bfy)\in\bbR\times H^{-1/2}(\partial D_\bfy)\setminus\{0\}$ such that \cref{eq:sloevp} holds is an equivalent reformulation of \cref{eq:LaplaceEVP} into a parameter dependent nonlinear eigenvalue problem. \cite[Lemma 2.4ff.]{SU2012} shows that, for each $\bfy$, the maps
\begin{align*}
\tilde{\calV}_\bfy\colon\bbC&\to\calL(H^{-1/2}(\partial D_\bfy),H^1(D_\bfy))\\
\kappa&\mapsto\tilde{\calV}_\bfy(\kappa)
\end{align*}
and
\begin{align*}
\calV_\bfy\colon\bbC&\to\calL(H^{-1/2}(\partial D_\bfy),H^{1/2}(\partial D_\bfy))\\
\kappa&\mapsto\calV_\bfy(\kappa)
\end{align*}
are holomorphic. For the reformulation of the eigenvalue problem to operators from and to $H^{1/2}(\partial D_\bfy)$ using Riesz isomorphisms we refer to \cite{SU2012}. Moreover, the results from \cite{DH2024} imply that the parameter-to-operator map $\bfy\mapsto(\kappa\mapsto\bfV_\bfy(\kappa))$ is holomorphic, justifying the use of quasi-Monte Carlo techniques.

\subsection{Modelling of the domain deformations}\label{subsec:modelling-domain}
We consider the family of param\-e\-trized domain boundaries $\{\partial D_\bfy\}_{\bfy\in Y}$ as a deformation of the unit cube $D_{\text{ref}}=(0,1)^3$ with boundary $\partial D_{\text{ref}}$. That is, we consider $\{\partial D_\bfy\}_{\bfy\in Y}$ to be given pointwise by
\[
\partial D_\bfy
=
\partial D_{\text{ref}}+\bfV(D_{\text{ref}},\bfy),
\]
where $\bfV\colon \partial D_{\text{ref}}\times Y\to\bbR^3$ is a random vector field. We set $Y=[-\frac{1}{3}, \frac{1}{3}]^{51}$ and prescribe the vector field to be given by a Karhunen-Lo\`eve expansion
\[
\bfV(\bfx,\bfy)
=
\sum_{j=1}^{51} \sqrt{\mu_j}\psi_j(\bfx) y_j,
\qquad
\bfx\in \partial D_{\text{ref}},~\bfy\in Y,
\]
where $\{(\mu_j,\psi_j)\}_{j=1}^{51}$ are the eigenpairs to the $51$ largest eigenvalues of the (compact) covariance operator $\calC\colon\big[L^2(\partial D_{\text{ref}})\big]^3\to \big[L^2(\partial D_{\text{ref}})\big]^3$:
\[
(\calC\varphi)(\bfx)
=
\int_{\partial D_{\text{ref}}}\operatorname{Cov}[\bfV](\bfx,\bfz)\varphi(\bfz)\dd\sigma_\bfz,
\]
where $\operatorname{Cov}[\bfV]\colon\partial D_{\text{ref}}\times\partial D_{\text{ref}}\to\bbR^{3\times 3}$ is the covariance function of the random vector field $\bfV$, which we prescribe as
\[
\operatorname{Cov}[\bfV](\bfx,\bfz)
=
\begin{bmatrix}
\frac{1}{10}e^{-\frac{\|\bfx-\bfz\|_2^2}{2\cdot 16^2}} & 0 & \frac{1}{10^6}e^{-\frac{\|\bfx-\bfz\|_2^2}{2\cdot 64^2}} \\
0 & \frac{1}{10}e^{-\frac{\|\bfx-\bfz\|_2^2}{2\cdot 4^2}} & 0 \\
\frac{1}{10^6}e^{-\frac{\|\bfx-\bfz\|_2^2}{2\cdot 64^2}} & 0 & \frac{1}{5}e^{-\frac{\|\bfx-\bfz\|_2^2}{2}}
\end{bmatrix}.
\]
For a more exhaustive discussion of the domain modeling part we refer to \cite{DHJM2022}.

The eigenvalues of \eqref{eq:LaplaceEVP} on the reference domain $D_{\text{ref}}$ are given by 
\[
    \lambda_k = \pi^2(k_1^2+k_2^2+k_3^2)
\]
with associated eigenfunctions
\[
    u_k(\bfx) = \sin(\pi k_1 x_1)\sin(\pi k_2 x_2)\sin(\pi k_3 x_3)
\]
for all $k = (k_1,k_2,k_3)\in \bbN^3$.

\subsection{Computational setup}

For the modelling of the domain and the computation of the discretizations of $\calV_\bfy(\kappa_\bfy)$, we employ the \texttt{C++} library \texttt{bembel} \cite{Bembel2020}. To solve the discretized version of \eqref{eq:sloevp} using the contour integral method, we choose a circle with center $\mu = \frac{(\sqrt{3}+\sqrt{6})\pi}{2}$ and radius $r = 2$.
In the reference case, it contains $\kappa_{(1,1,1)} = \sqrt{3}\pi, \kappa_{(1,1,2)} = \kappa_{(1,2,1)} = \kappa_{(2,1,1)} = \sqrt{6}\pi$, i.e., $m=4$ eigenvalues counting algebraic multiplicities.
As in the previous example, the second eigenvalue, which is triple, does not remain a repeated eigenvalue over the entire parameter space but the trajectories bifurcate generating three distinct, simple eigenvalues for parameters $\bfy$ different from zero.

Due to the high asymptotic convergence rate of $h^{2p+3}$ of the boundary element method for the approximation of the eigenvalues, the number of samples for the quasi-Monte Carlo quadrature has to be adapted for each level according to $\sim h^{-2p-3}$.
This yields the number of quadrature points as shown in \cref{tab:number-of-samples_BEM}.

\begin{table}[htb]
    \centering
\begin{tabular}{c c c c c c}
    \hline
            & $\ell=2$    & $\ell=3$ & $\ell=4$ & $\ell = 5$ & $\ell = 6$\\\hline
    $p = 0$ & 409\,600    & 50\,120  & 6\,400   & 800        & 100 \\\hline
    $p = 1$ & 3\,276\,800 & 102\,400 & 3\,200   & 100        & -- \\\hline
    \end{tabular}
    \medskip
    \caption{Number of QMC samples per level for the nonlinear formulation of the Laplacian eigenvalue problem.}
    \label{tab:number-of-samples_BEM}
\end{table}

The computations were carried out in parallel on eight nodes of the aforementioned cluster, making 512 cores in total.

\subsection{Numerical results}
\begin{figure}
    \centering
    \begin{tikzpicture}[scale=0.72]
        \begin{groupplot}[
            group style={
                group size=3 by 1,
                horizontal sep=1.5cm,
            },
            width=6cm,
            height=5cm,
            ymode=log,
            xlabel={$L$},
            legend style={font=\footnotesize}
        ]
        \nextgroupplot[title={Eigenvalues (Eulerian)}, xtick={0,1,2,3,4}, xticklabels={0,1,2,3,4},legend pos=north east]
        \addplot[color=red,mark=x]table[x=L,y=eveuler]{figures/bem_p0.data};\addlegendentry{$p=0$}
        \addplot[color=red, mark=square]table[x=L,y=eveuler]{figures/bem_p1.data};\addlegendentry{$p=1$}
        \addplot[domain=0:4.2, samples=200, dashed]{0.2*2^(-3*x)}; \addlegendentry{$h^{2p+3}$}
        \addplot[domain=0:3.2, samples=200, dashed]{0.02*2^(-5*x)};

        \nextgroupplot[title={Eigenvalues (Lagrangian)}, xtick={0,1,2,3,4}, xticklabels={0,1,2,3,4},legend pos=north east]
        \addplot[color=red,mark=x]table[x=L,y=evlagrange]{figures/bem_p0.data};\addlegendentry{$p=0$}
        \addplot[color=red, mark=square]table[x=L,y=evlagrange]{figures/bem_p1.data};\addlegendentry{$p=1$}
        \addplot[domain=0:4.2, samples=200, dashed]{0.2*2^(-3*x)}; \addlegendentry{$h^{2p+3}$}
        \addplot[domain=0:3.2, samples=200, dashed]{0.02*2^(-5*x)};
        
        \nextgroupplot[title={Eigenspaces (Lagrangian)}, xtick={0,1,2,3,4}, xticklabels={0,1,2,3,4}, legend pos=south west]
        \addplot[color=red,mark=x]table[x=L,y=eslagrange]{figures/bem_p0.data}; \addlegendentry{$p=0$}
        \addplot[color=red, mark=square]table[x=L,y=eslagrange]{figures/bem_p1.data}; \addlegendentry{$p=1$}
        \addplot[domain=0:4.2, samples=200, dashed]{1.1*2^(-1*x)}; \addlegendentry{$h^{p+1}$}
        \addplot[domain=0:3.2, samples=200, dashed]{0.2*2^(-2*x)};

        \end{groupplot}
    \end{tikzpicture}
    \caption{Convergence of multilevel QMC for the expectation of the Eulerian and Lagrangian parameter-to-solution maps.}
    \label{fig:plots}
\end{figure}
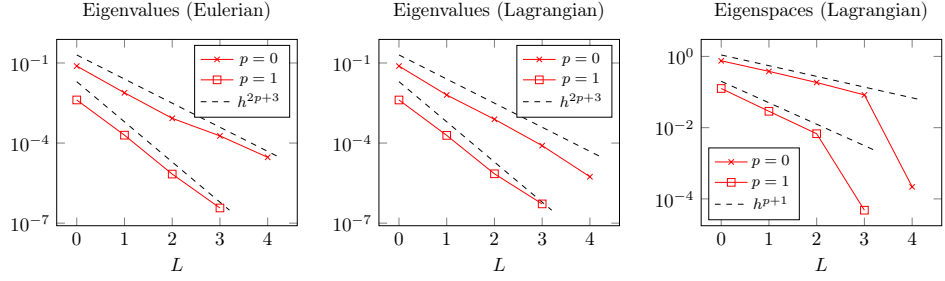

In order to measure the approximation errors, we proceed analogously to \cref{subsec:toyproblem-numerical-results}.
\Cref{fig:plots} shows the convergence of the parameter-to-solution maps computed with the Kronecker sequence towards a reference solution computed using the Halton sequence.
As we are using the $L^2$-norm to measure the difference between the eigenspaces, instead of the $h^{p+3/2}$ that we would obtain for the $H^{-1/2}$-norm we expect a rate of convergence of $h^{p+1}$, which is the observed convergence in \Cref{fig:plots}.

\begin{figure}
    \centering
    \includegraphics[width=0.7\textwidth]{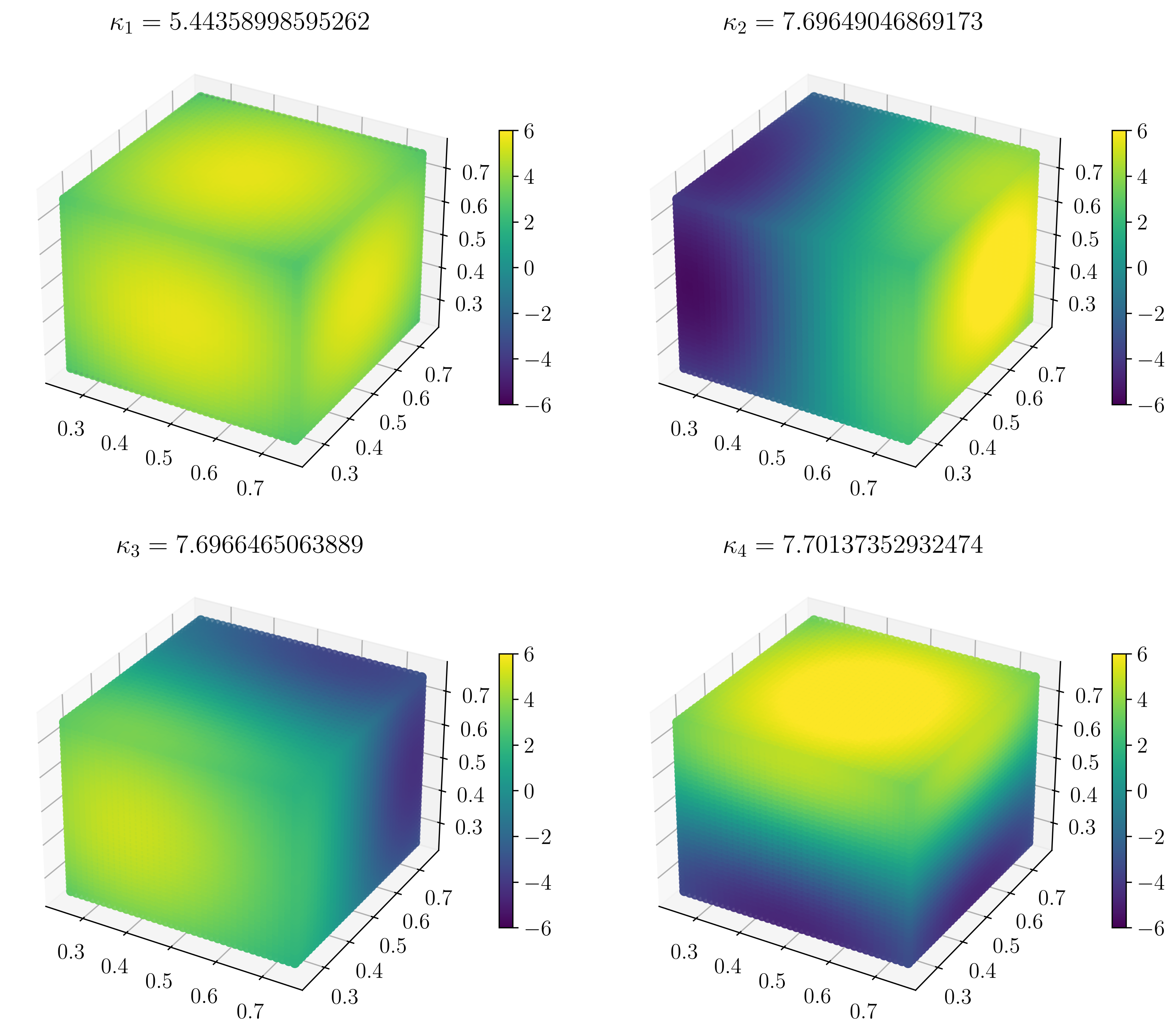}
    \caption{Visualization of the expected eigenfunctions (Lagrangian para\-meter-to-eigenspace map) for the 3D Laplacian eigenvalue problem evaluated on the surface of the cube $[0.25,0.75]^3 \subset D_{\text{ref}}$.}
    \label{fig:visual-3d-exp-eigfct}
\end{figure}

\Cref{fig:visual-3d-exp-eigfct} illustrates the expected eigenfunctions obtained from the Lagrangian parameter-to-eigenspace map.

\section{Conclusion}\label{sec:conclusion}
Statistical quantities of interest of random eigenvalue problems are difficult to define and compute in general, since crossings and bifurcations make parameter-to-solution maps not well defined in general. In this article, we proposed parameter-to-solution maps  for random nonlinear eigenvalue problems from an Eulerian and from a Lagrangian perspective. Both of these maps have straightforward theoretical interpretations and can be easily implemented algorithmically by modification of Beyn's contour integral method. Using these maps and algorithic implementation, we proposed an algorithmic approach for multilevel Monte Carlo estimators which matching or tracking of eigenvalues is not required. Prolongation operators are only required when the corresponding eigenspaces are of interest. The algorithmic approach is adaptable to other quadrature methods such as quasi-Monte Carlo methods or sparse grids. A limitation of the multilevel approach is that the convergence theory for $m>1$, i.e., for clusters of eigenvalues and their generalized eigenspaces, seems not to be available and needs to be presumed from the $m=1$ case. Numerical experiments for the random diffusion eigenvalue problem and a boundary integral formulation of the Laplacian eigenvalue problem demonstrate this presumption empirically. Other boundary integral formulations such as for Maxwell \cite{Ung2021} or fluid-structure interaction \cite{EKSU2012} could be considered likewise.

\bibliographystyle{plain}
\bibliography{references}

\end{document}